\documentclass[12pt]{article}
\usepackage{amssymb}
\newcommand{\sect}[1]{\setcounter{equation}{0}\section{#1}}

\def \d2dots{\mathinner{\mkern1mu\raise1pt\vbox{\kern7pt\hbox{.}}\mkern2mu
\raise4pt\hbox{.}\mkern2mu\raise7pt\hbox{.}\mkern1mu}}

\newcommand{\be}{\begin{equation}}
\newcommand{\ee}{\end{equation}}
\newcommand{\bea}{\begin{eqnarray}}
\newcommand{\eea}{\end{eqnarray}}
\newcommand{\beano}{\begin{eqnarray*}}
\newcommand{\eeano}{\end{eqnarray*}}
\newcommand{\nonu}{\nonumber \\}
\newcommand{\vertequiv}{\rule{0.0642ex}{1.42ex}\,\rule{0.0642ex}{1.42ex}\,
                      \rule{0.0642ex}{1.42ex}}
\newcommand{\eps}{\epsilon}

\newcommand{\vph}{\varphi}

\newcommand{\ca}{\mbox{$\cal{A}$}}

\newcommand{\cd}{\mbox{$\cal{D}$}}

\newcommand{{\cg}}{\mbox{$\cal{G}$}}
\newcommand{\ch}{\mbox{$\cal{H}$}}
\newcommand{\ci}{\mbox{${\cal I}$}}

\newcommand{\cl}{\mbox{$\cal{L}$}}

\newcommand{\ct}{\mbox{$\cal{T}$}}
\newcommand{\cu}{\mbox{${\cal U}$}}
\newcommand{\cw}{\mbox{$\cal{W}$}}

\def\opluslim{\mathop{\oplus}\limits}
\def\otimelim{\mathop{\otimes}\limits}

\newcommand{\wh}[1]{\widehat{#1}}

\newcommand{\half}{\frac{1}{2}}

\newtheorem{prop}{Property}[section]
\newtheorem{coro}[prop]{Corollary}
\newtheorem{defi}[prop]{Definition}
\newtheorem{theor}[prop]{Theorem}

\newtheorem{lem}[prop]{Lemma}
\newcommand{\prf}{\underline{Proof:}\ }
\newcommand{\finprf}{\null \hfill {\rule{5pt}{5pt}}\\[2.1ex]\indent}
\newcommand{\ie}{{\it i.e.}\ }
\newcommand{\1}{\mbox{\hspace{.0em}1\hspace{-.24em}I}}

\newcommand{\CC}{\mbox{${\mathbb C}$}}
\newcommand{\BB}{\mbox{${\mathbb B}$}}
\newcommand{\II}{\mbox{${\mathbb I}$}}
\newcommand{\JJ}{\mbox{${\mathbb J}$}}
\newcommand{\JJc}{\mbox{${\mathbb J}\vert_{const.}$}}
\newcommand{\JJgf}{\mbox{${\mathbb J}\vert_{g.f.}$}}

\newcommand{\NN}{\mbox{${\mathbb N}$}}
\newcommand{\WW}{\mbox{${\mathbb W}$}}

\newcommand{\CMP}[1]{Commun.\ Math.\ Phys.\ {\bf #1}}

\newcommand{\ynp}[1]{\mbox{$Y(N)_#1$}}
\newcommand{\synp}[2]{\mbox{$S_{#1}Y(N)_{#2}$}}

\begin{document}
\renewcommand{\thefootnote}{\fnsymbol{footnote}}
\newpage
\pagestyle{empty}
\setcounter{page}{0}

%%%%%%%%%%%%%%%%%%%%%%%%%%%%%%%%
%%%%%  DEBUT PAGE DE TITRE  %%%%%%
%%%%%%%%%%%%%%%%%%%%%%%%%%%%%%%%
%%%%%%%%%%%% LOGO LAPTH - DEBUT  %%%%
%%%%%%%%%%%%%%%%%%%%%%%%%%%%%%%%

\newcommand{\LAP}{LAPTH}
\def\logo{{\bf {\huge LAPTH}}}

\centerline{\logo}

\vspace {.3cm}

\centerline{{\bf{\it\Large 
Laboratoire d'Annecy-le-Vieux de Physique Th\'eorique}}}

\centerline{\rule{12cm}{.42mm}}
%%%%%%%%%%%%%%%%%%%%%%%%%%%%%%%%%%%%%%%
%%%%%%%%%%%%%%%% LOGO LAPTH  - FIN %%%%%%%%%%
%%%%%%%%%%%%%%%%%%%%%%%%%%%%%%%%%%%%%%%

\vfill

\begin{center}
{\LARGE\bf RTT presentation of finite \cw-algebras}\\[2.1em]

{\large
C. Briot\footnote{briot@lapp.in2p3.fr} and 
E. Ragoucy\footnote{ragoucy@lapp.in2p3.fr}}

\null

LAPTH\footnote{UMR 5108 du CNRS associ\'ee \`a l'Universit\'e de 
Savoie.}, Chemin de Bellevue, BP 110, F-74941 Annecy-le-Vieux
  cedex, France
\end{center}

\vfill

\begin{abstract}
    We construct a wide class of finite \cw-algebras as truncations of 
    Yangians. These truncations correspond to algebra homomorphisms
    and allow to construct the \cw-algebras as exchange algebras, the 
    R-matrix being the Yangian's one. 
    
    As an application, we classify all irreducible finite 
    dimensional representations of these \cw-algebras and determine 
    their center.
\end{abstract}

\vfill
\rightline{April 2000}
\rightline{LAPTH-792/00}
\rightline{\tt math.QA/0005111}
\newpage
\pagestyle{plain}
\setcounter{page}{1}

%%%%%%%%%%%%%%%%%%%%%%%%%%%%%%%%
%%%%%  FIN PAGE DE TITRE  %%%%%%
%%%%%%%%%%%%%%%%%%%%%%%%%%%%%%%%

%%%%%%%%%%%%%%%%%%%%%%%%%%%%%%%%
%%%%%  HEADINGS POUR DRAFT  %%%%%%

%\markright{\today\dotfill DRAFT\dotfill }
%\pagestyle{myheadings}

%%%%%%%%%%%%%%%%%%%%%%%%%%%%%%%%

\setcounter{footnote}{0}
\tableofcontents
\newpage
\sect{Introduction}

It has already been proven \cite{RS} that there exists an algebra 
homomorphism between Yangian based on $sl(N)$ and finite 
$\cw(sl(Np),N.sl(p))$-algebras. 
Such a connection plays a role in the study of physical models: for 
instance, in the case of the $N$-vectorial non-linear Schr\"odinger equation on the 
real line, the full symmetry is the Yangian $Y(gl(N))$, but the space of states 
with particle number less than $p$ is a representation of the 
$\cw(gl(Np),p.sl(N))$ algebra \cite{schro}.

The connection between Yangians and finite $\cw(sl(Np),N.sl(p))$-algebras 
was proven in the Drinfeld presentation \cite{Drinfeld} of the 
Yangian. Since the homomorphism is (obviously) not an isomorphism, it 
does not allow to carry the Yangian R-matrix "down" to the finite
 \cw-algebras. In this paper, we prove the correspondence in the 
"RTT" presentation \cite{YRTT}of the Yangian. The mentioned finite 
\cw-algebras\footnote{More precisely, it is the $\cw(gl(Np),N.sl(p))$ algebras 
which are concerned, we will come back on this slight difference later 
on.} appear to be 
"truncation" of the Yangian, \ie the resulting coset when modding out 
the Yangian "high level" generators. These truncated Yangians were 
already introduced in \cite{Chered} under the name of Yangian of 
level $p$.
Thanks to this presentation, we 
can deduce a R-matrix for the \cw-algebras under consideration, as 
well as the complete classification of the finite-dimensional 
irreducible
representations of these algebras. We also show that the Hopf 
structure of the Yangian  cannot be carried by the homomorphism: 
although this is 
not a "no-go theorem" for \cw-algebras to be Hopf algebras, 
 it severely constrains the possibilities to get this structure.

\null

To prove our result, we need to combine three notions: Yangians, 
\cw-algebras and cohomology. We have tried to be self-contained, and 
as such, we need to recall known results for these different fields: it 
is done in section \ref{yang} for Yangians, in section \ref{walg} for
\cw-algebras, and in the appendix \ref{coho} for cohomology. We 
collect our results in section \ref{reslt} and then present
applications in section \ref{appl}. We conclude with a \ldots 
conclusion, where possible generalizations and applications of our 
results are presented (section \ref{conclu}). Some calculations about $gl(Np)$ 
algebras  are collected in the appendix \ref{notations}.

\sect{Yangians\label{yang}}
Yangians can be seen as deformations of loop algebras (based on a 
simple Lie algebra) and associated to a rational solution to the 
Yang-Baxter equation. They have been extensively studied, and we 
refer to \cite{Drinfeld,Chari,MNO} and references therein for more details. 
We will here focus on Yangians based on $gl(N)$, and recall the basic 
properties below.

\subsection{The Yangian $Y(gl(N))$}
There is essentially two presentations of $Y(gl(N))$: one based on 
generators and relations \cite{Drinfeld} (Serre-Chevalley-type presentation), 
and the second (closer to integrable systems methods) using the $R$-matrix 
approach \cite{YRTT} (see also \cite{MNO,Chari} and ref. therein). 
We use here the last one. The generators of the Yangian 
are gathered in a single matrix:
\begin{equation}
    T(u)=\sum_{n=0}^\infty\sum_{i,j=1}^N u^{-n} T^{ij}_{n} E_{ij}=
    \sum_{n=0}^\infty u^{-n} T_{n} 
    =\sum_{i,j=1}^N  T^{ij}(u) E_{ij} 
    \mbox{ with } T^{ij}_{0}=\delta^{ij}
    \label{eq:defT}
\end{equation}
where $u$ is a spectral parameter and $i,j$ indices in the fundamental 
of $gl(N)$. $E_{ij}$ is the usual $N\,$x$\, N$ matrix with 1 at position 
$(i,j)$. The algebraic structure is encoded in the relation
\begin{equation}
    R(u-v)T_{1}(u)T_{2}(v) = T_{2}(v)T_{1}(u)R(u-v)
    \label{eq:RTT}
\end{equation}
with $R(x)=1\otimes 1-\frac{1}{x}P_{12}$ and $P_{12}$ is the flip 
operator ($P_{12}=\sum_{i,j=1}^N  E_{ij}\otimes E_{ji}$ in 
representations). 
The commutation relations read in components
\begin{equation}
    [T^{ij}_{m}{,} T^{kl}_{n}] = \sum_{r=0}^{\mbox{\scriptsize min}(m,n)-1} (
    T^{kj}_{r}T^{il}_{m+n-r-1}-T^{kj}_{m+n-r-1}T^{il}_{r})
    \label{comTT}
\end{equation}
Note that in (\ref{comTT}), all the couples $(r,\ s)$, where $s=m\!+\!n\!-\!1\!-\!r$, satisfy 
$s<$min$(m,n)$ and $r\geq$max$(m,n)$.

It is known that the  Yangian $Y(N)$ is a deformation of a loop algebra 
based on $gl(N)$. The parameter $\hbar$ can be recovered by 
multiplying the generators by an appropriate power of $\hbar$:
\be
T^{ij}_{n}\ \rightarrow\ \hbar^{n-1} T^{ij}_{n}
\ee
Then, the relations (\ref{comTT}) can be rewritten as
\be
    [T^{ij}_{m}{,} T^{kl}_{n}] = 
    \delta^{kj} T^{il}_{m+n-1}-\delta^{il} T^{kj}_{m+n-1} +o(\hbar)
\ee
which shows that $Y(N)$ is a deformation  of a loop algebra 
(restricted to its positive modes).
It can be proven that as soon as $\hbar\neq0$, all the Hopf algebras 
$Y_{\hbar}(N)$ are isomorphic. 

The Hopf structure is given by
\be
    \Delta(T(u))  =  T(u)\otimes T(u)
    \label{eq:deltaT}  \ \ ; \ \ 
    \epsilon(T(u))  =  1   
     \ \ ; \ \ 
    S(T(u)) =  -T(u)
\ee
or in components:
\be
    \Delta(T^{ij}_{m})  = \sum_{k=1}^N\sum_{r=0}^mT^{ik}_{r}\otimes 
    T^{kj}_{m-r}
    \label{comp:deltaT}  \ \ ; \ \ 
    \epsilon(T^{ij}_{m})  =  \delta_{m,0}\delta_{i,j}   
     \ \ ; \ \ 
    S(T^{ij}_{m})  =  -T^{ij}_{m}
\ee
For briefness, we will denote $Y(N)\equiv Y(gl(N))$.

\subsection{Center of $Y(N)$ and associated Hopf 
subalgebras.\label{centY}}

The center $\cd=\cu(d_{i},\ i\in\NN)$ of $Y(N)$ is generated by the quantum determinant:
\begin{equation}
    \mbox{q-det}(T(u))\equiv\sum_{\sigma\in{\cal 
\Sigma}}\mbox{sgn}(\sigma)T_{\sigma(1)1}(u)T_{\sigma(2)2}(u-1)\cdots 
T_{\sigma(N)N}(u-N+1)=1+\sum_{n=1}^\infty u^{-n} d_{n} 
%\mbox{ with } d_{0}=1
    \label{eq:qdet}
\end{equation}
The Hopf algebra $Y(sl(N))$ is the quotient of $Y(N)$ by the 
relation q-det$T=1$, \ie $Y(sl(N))\sim Y(N)/\cd$.

We introduce 
$$
{\cal D}_{r}=\cu(\{d_{1},d_{2},\dots,d_{r}\})
$$
It is not 
difficult to show that for any value of $r$, $\cd_{r}$ is a Hopf ideal 
of $Y(N)$. It is obviously an algebra ideal (because it lies in
the center of the Yangian), and from (\ref{comp:deltaT}), one shows 
that 
\begin{equation}
    \Delta({\cal D}_{r})\subset{\cal D}_{r}\otimes {\cal 
    D}_{r}\ \Rightarrow\ 
    \Delta({\cal D}_{r})\ \subset\ {\cal D}_{r}\otimes Y(N)\oplus 
    Y(N)\otimes {\cal D}_{r}
\end{equation}
 hence $\cd_{r}$ is a coideal. Consequently, the coset 
$Y(N)/{\cal D}_{r}$ is also a Hopf algebra.
\begin{equation}
    \synp{r}{}=Y(N)/{\cal D}_{r}\mbox{ and }Y(N)\sim 
    \synp{r}{}\otimes {\cal D}_{r}
\end{equation}
This allows us to construct a series of Hopf subalgebras:
\beano
&&\synp{r}{}=Y(N)/{\cal D}_{r}\mbox{ and }Y(N)\sim 
    \synp{r}{}\otimes {\cal D}_{r}\ \ \forall r\\
&&Y(N)\equiv\synp{0}{} \  \supset \  \synp{1}{} \  \supset \  
        \cdots \  \supset \  \synp{r}{} \  \cdots \  \supset \  
        \synp{\infty}{}\equiv Y(sl(N))  
\eeano
where $Y(sl(N))$ is the only one which possesses a trivial center. 
The intermediate subalgebras will be of some use in the following.

\subsection{Evaluation representations \label{sec.eval}}
The finite dimensional irreducible representations of $Y(N)$ have been 
classified \cite{Dreval,tara}, see also \cite{Chari,Mol} for more 
details. 
It uses the notion of evaluation representations\cite{KRS,CPeval}:
\begin{defi}\label{def.eval}
    {\bf Evaluation representations}\\    
    An evaluation representation $ev_{\pi}$ is a morphism from the 
    Yangian $Y(gl(N))$ to a highest weight irreducible
    representation $\pi$ of $gl(N)$. 
    The morphism is given by
    \be
    ev_{\pi}(T^{ij}_{(1)})=\pi(T^{ij}_{(1)}) 
%    \ \mbox{ ; }\ ev_{\pi}(T^{ij}_{(1)})=A\,\pi(T^{ij}_{(0)})
    \ \mbox{ and }\ ev_{\pi}(T^{ij}_{(n)})=0,\ n>1
    \ee
    where we have identified the generators $T^{ij}_{(1)}$ with 
    $gl(N)$ elements.
\end{defi}
The evaluation representations form a very simple class of 
representations, since only one kind of Yangian generators
is non-trivially represented. They are sufficient to get all 
finite-dimensional irreducible  representations, through the tensor products of such 
representations:
\begin{defi}\label{tens.eval}
    {\bf Tensor product of evaluation representations}\\    
    Let  $\{ev_{\pi_{i}}\}_{i=1,..,n}$ be a set of evaluation 
    representations. The tensor product of these 
    $n$ representations  
    $ev_{\vec{\pi}}=ev_{\pi_{1}}\otimes..\otimes ev_{\pi_{n}}$ is a morphism from the 
    Yangian $Y(gl(N))$ to the tensor 
    product of $gl(N)$ 
    representations $\vec{\pi}=\otimes_{i}\pi_{i}$  given by
    \be
       ev_{\vec{\pi}}(T^{ij}_{(r)})=\opluslim_{r_{1}+r_{2}+..+r_{n}=r}\left(
     \otimelim_{k=1}^n\ ev_{\pi_{k}}(T^{ij}_{(r_{k})})\right)
    \ee
    It satisfies:
    \be
    ev_{\vec{\pi}}(T^{ij}_{(r)})\neq0\ \mbox{ if and only if }\ r\leq n
    \label{tensW}
    \ee
\end{defi}
Note that this definition follows from the Yangian coproduct (\ref{eq:deltaT}). 
Tensor product of
evaluation representations play an important role in the 
classification of finite dimensional irreducible 
representations of Yangians. This is reflected in the following 
theorems and corollary (proved in \cite{Dreval}, see also \cite{Mol,plus,tara} 
for more details).

\null 

{\bf Theorem:} {\em Any finite dimensional irreducible 
representation of $Y(N)$ is highest weight and contains (up to 
multiplication by a scalar) a unique highest weight vector.}

\null 

By highest weight vector, we mean a vector $\eta$ (in the 
representation) such that 
\[
\begin{array}{ll}
t^{ij}(u)\eta = 0 &1\leq i<j\leq N \\
t^{ii}(u)\eta = \lambda^i(u)\, \eta &  1\leq i\leq N 
\end{array}
\]
where $\lambda^i(u)=1+\sum_{r>0}\lambda_{(r)}^i u^{-r}$, with
$\lambda^i_{(r)}\in\CC$,  and $t^{ij}(u)$ 
represents $T^{ij}(u)$. 
As usual, $\lambda(u)=(\lambda^1(u),\ldots,\lambda^N(u))$ is called 
the weight of the representation.

\null 

{\bf Theorem:} {\em An irreducible highest weight
representation of $Y(N)$ of weight $\lambda(u)$ is finite dimensional 
if and only if there exist $(N-1)$ monic polynomials $P_{i}(u)$ such that 
\[
\frac{\lambda^i(u)}{\lambda^{i+1}(u)}=\frac{P_{i}(u+1)}{P_{i}(u)}
\]

In that case, the representation is isomorphic to the subquotient of 
the tensor product of $m=\sum_{i}m_{i}$ evaluation representations, 
where $m_{i}$ is the degree of $P_{i}(u)$.}

\null 

By monic polynomials we mean a polynomial of the form
\[
P_{i}(u)=\prod_{k=1}^{m_{i}} (u-\gamma_{k})\ \mbox{ with }\ 
\gamma_{k}\in\CC
\]
By subquotient, we mean the irreducible part of the highest weight submodule 
of the mentioned tensor 
product. More precisely, in the tensor product of  
evaluation
representations (which are by definition highest weight 
representations), one considers the submodule generated by the tensor 
product 
of the highest weight vectors, and quotients it by all (sub)singular vectors 
which may appear. 

Note that although generically
the tensor product is  irreducible (\ie is equal to the mentioned submodule and has no 
singular vector),  it is 
only for $Y(2)$ that it is {\em always} irreducible (see 
counter-example for $Y(3)$ in \cite{Mol}).

A simpler characterization of the finite dimensional irreducible 
representations is given by the following corollary

\null

{\bf Corollary} {\em The irreducible finite dimensional 
representations of $Y(N)$ are in one-to-one correspondence with the 
families $\{P_{1}(u),\ldots,P_{N-1}(u),\rho(u)\}$ where $P_{i}$ are monic 
polynomials and $\rho(u)=1+\sum_{n>0}d_{n}u^{-n}$ encodes the values of 
the central elements.}

\subsection{Truncated Yangians}
The notion of truncated Yangians has been already introduced in 
\cite{Chered} (although not named truncated, but Yangian of level $p$) as a tool in  
 representation theory. They were also studied in \cite{yangTr}. 
We now introduce the left ideal generated by $\ct_{p}=\cu(\{T^{ij}_{n}, 
n>p\})$ :
$$
\ci_{p}=Y(N)\cdot\ct_{p}
$$
 and  the coset (truncation of the Yangian at order $p$)
 \be
 \ynp{p}=Y(N)/\ci_{p}
 \ee

{\prop The truncated Yangian $\ynp{p}$ is an 
algebra ($\forall N\in\NN$, $\forall p\in\NN$). $\Delta$ is not a 
morphism of this algebra (for the structure induced by 
$Y(N)$).}

\null

\prf
We prove the Lie algebra structure of $\ynp{p}$ by showing that $\ct_{p}$ is a
bilateral  ideal, \ie that we have $Y(N)\cdot \ci_{p}\subset 
\ci_{p}$. In fact, we will show a more stronger property, that is
\be
\left[ Y(N)\, , \, \ct_{p}\right]\subset Y(N)\cdot \ct_{p}
\ \mbox{ and }\ 
\left[ Y(N)\, , \, \ct_{p}\right]\subset \ct_{p}\cdot Y(N) 
\label{eq:ideal}
\ee 
We make the calculation for the first inclusion, 
the proof for the other inclusion being identical.
Indeed, the relation (\ref{comTT}) shows that 
$[T^{ij}_{m},T^{kl}_{n}]$ (for $n>p$) is the sum of two terms, 
the first being in $Y(N)\cdot \ct_{p}$, the second belonging to 
$\ct_{p}\cdot Y(N)$. Focusing on the latter, one rewrites it as
\be
\begin{array}{l}
\displaystyle\hspace{-2.4ex}
\sum_{r=0}^{\mbox{\scriptsize min}-1} T^{kj}_{m+n-1-r} T^{il}_{r} \ =\ 
\sum_{r=0}^{\mbox{\scriptsize min}-1}\left(T^{il}_{r} T^{kj}_{m+n-1-r} +
\sum_{s=0}^{r-1}\left(T^{ij}_{s}T^{kl}_{m+n-2-s}-T^{ij}_{m+n-2-s}T^{kl}_{s}\right)
\right) \\
\displaystyle\hspace{1.8em} =\  
\sum_{r=0}^{\mbox{\scriptsize min}-1}T^{il}_{r} T^{kj}_{m+n-1-r} +
\sum_{s=0}^{\mbox{\scriptsize min}-2}(\mbox{min}-s-1)\left(
T^{ij}_{s}T^{kl}_{m+n-2-s}-T^{ij}_{m+n-2-s}T^{kl}_{s}\right) 
\end{array}
\label{TT=TT}
\ee
where min stands for min$(m,n)$. In (\ref{TT=TT}), 
only the last term belongs to $\ct_{p}\cdot Y(N)$, 
with a summation which has one term less than the previous one:
we can thus proceed recursively in a finite number of steps. The final result 
is an element of $Y(N)\cdot \ct_{p}$.

\null

As far as Hopf structure is concerned, the calculation
$$
\Delta(T^{ij}_{p+1})=T^{ij}_{p+1}\otimes1+1 \otimes T^{ij}_{p+1}
+\sum_{n=1}^{p}T^{ik}_{n}\otimes T^{kj}_{p+1-n}
$$
shows that $\ci_{p}$ is not a coideal, since we have
$$
\Delta(\ci_{p})\ \not\subset\ Y(N)\otimes \ci_{p}\oplus \ci_{p}\otimes 
Y(N)
$$
Moreover, $\Delta$ is not an algebra morphism 
anymore, since for instance
\be
\Delta\left(\left[T^{ij}_{p},T^{kl}_{2}\right]\right)-
\left[\Delta(T^{ij}_{p}),\Delta(T^{kl}_{2})\right]=
\sum_{s+t=p}(T^{il}_{s+1}\otimes T^{kj}_{t}-T^{il}_{s}\otimes 
T^{kj}_{t+1})
\neq 0
\ee
\finprf
Finally, we note that each $\ynp{p}$ is a deformation of a truncated 
loop algebra based on $gl(N)$. By truncated loop algebra, we mean the 
quotient  of a usual $gl(N)$ loop algebra (of generators $t^{ij}_{n}$) by the 
relations $t^{ij}_{n}=0$ for $n<0$ and $n>p$.
The construction is the same as  for the 
complete Yangian.

\subsection{Poisson Yangians}
In the following we will deal with a Poisson version of the Yangian, 
where the commutator is replaced by Poisson bracket. It corresponds 
to the usual classical limit of quantum groups. One sets
\be
T(u)= L(u)\ ;\ R_{12}(x)=\1 +\hbar\, 
r_{12}(x)+o(\hbar)
\ ;\ [\ ,\ ]=\hbar\{\ ,\ \}+o(\hbar)
\ee
The relation (\ref{eq:RTT}) is then expanded as a series in $\hbar$, the 
first non-trivial term being the $\hbar^2$ coefficient. This new 
relation is the defining 
relation for the Poisson Yangian and reads:
\begin{equation}
    \left\{L(u) \stackrel{{\otimes}}{,} L(v)\right\} = [r_{12}(u-v), 
    L(u)\otimes L(v)]
    \mbox{ with } r_{12}(x)=\frac{1}{x} P_{12}
\end{equation}
where $\{L(u)\! \stackrel{{\otimes}}{,}\! L(v)\}$ is a matrix of 
component $\{L^{ij}{,} L^{kl}\}$ in the basis $E_{ij}\otimes E_{kl}$. 
In components
\begin{equation}
    \{T^{ij}_{m} , T^{kl}_{n}\} = \sum_{r=0}^{\mbox{\scriptsize min}(m,n)-1} (
    T^{kj}_{r}T^{il}_{m+n-r-1}-T^{ij}_{m+n-r-1}T^{ij}_{r})
\end{equation}
Apart from the change from commutators to Poisson brackets (and the 
commutativity of the product), all the above algebraic properties still apply.

In particular, we can still define the truncated (Poisson) Yangian, 
with the same procedure as above.

\sect{\cw-algebras\label{walg}}
 Such algebras can be constructed by symplectic
reduction of finite dimensional Lie algebras in the same way the 
conformal (affine) \cw-algebras \cite{zam} arise as reduction of Kac-Moody (affine) 
Lie algebras \cite{ORaf}, 
 hence
the name finite \cw-algebras for the former \cite{Wfinie}. Some properties of such 
\cw-algebras have been developed \cite{Wrem}-\cite{Wanyon}. 
In particular, starting from a simple
Lie algebra $\cg$, a large class of \cw-algebras
 can be seen as the commutant, in a localization of the
enveloping algebra $\cu(\cg)$, of a \cg-subalgebra \cite{Wcom}. 
This feature  has already been exploited in various physical contexts
\cite{so42,Wanyon}. 
A remarkable fact is that the involved \cw-algebras 
are just of the type $\cw(sl(2n),n.sl(2))$, a subclass of the 
$\cw[gl(Np),N.sl(p)]$ algebras, in which we are interested here. 

We note $\cw_{p}(N)\equiv\cw[gl(Np),N.sl(p)]$. This algebra is 
defined as the Hamiltonian reduction of the enveloping algebra of 
$gl(Np)$ (see below). In general, the \cw-algebras are defined using 
semi-simple Lie algebras, but for $gl(m)$, we have the following property
$$
\cw[gl(m),\ch]\equiv \cw[sl(m)\oplus gl(1),\ch]\equiv \cw[sl(m),\ch]\oplus gl(1)
$$
which allows to extend the \cw-algebra to $gl(m)$. 

Note also that we are dealing with {\em finite} \cw-algebra, \ie the 
$gl(m)$ algebras we are speaking of are finite dimensional Lie 
algebras (not their affinization). 

We use the notations introduced in the appendix \ref{notations}.

\subsection{$\cw_{p}(N)$ as an Hamiltonian reduction\label{HamRed}}
Following the usual technic (see \cite{Wfinie} and \cite{Wanyon} for more 
details), we gather the generators of 
$gl(Np)$ in a $(Np)\times (Np)$ matrix:
\be
\JJ=\sum_{a,b=1}^{N}\sum_{j=0}^{p-1}\sum_{m=-j}^{j} 
J_{jm}^{ab} M^{jm}_{ab}
\ee
where $M^{jm}_{ab}$ are $(Np)\times (Np)$ matrices and $J_{jm}^{ab}$ 
are in the dual algebra of $gl(Np)$. They obey Poisson Brackets (PB) which mimic the 
commutation relations of $gl(Np)$:
\be
\{J_{ab}^{j,m}, J_{cd}^{\ell,n}\} = 
    \sum_{r=\vert j-\ell\vert}^{ j+\ell}\sum_{s=-r}^r \left(
    \rule{0ex}{2.4ex}
    \delta_{bc} <j,m;\ell,n\vert r,s> J_{ad}^{r,s} - 
    \delta_{ad} <\ell,n;j,m\vert r,s> J_{cb}^{r,s} \right)
    \label{PBgln}
  \ee
On the dual algebra, we introduce  first class constraints:
\be
\JJ\vert_{const.}=\eps_{-}+\sum_{a,b=1}^{N}\sum_{j=0}^{p-1}\sum_{m=0}^j
J_{jm}^{ab} M^{jm}_{ab}\equiv \eps_{-}+\BB
\ee
Explicitly, these constraints are imposed on the negative grade 
generators $J_{jm}^{ab}$, $m<0$, $\forall j,a,b$. They correspond to 
the vanishing of all these negative grade generators, but $J_{1,-1}^{00}$ 
which is set to 1. We will denote them generically by $\phi_{\mathbf 
x}$. 
Physically, these first class constraints generate gauge 
 transformations, an infinitesimal form of which is:
 \be
 \delta_{\lambda} J^{ab}_{jm} \sim \sum_{\mathbf x} \lambda_{\mathbf 
x}\ \{\phi_{\mathbf x}, J^{ab}_{jm}\}\label{gaugtr}
 \ee
 where the symbol $\sim$ means that one has to impose the constraints 
 once the PB has been computed. The interesting quantities are the 
 gauge invariant ones, and it can be shown that a way to construct a 
 basis for them is to choose a gauge fixing for \JJc. In the present 
 case, the gauge fixing is the highest weight gauge:
\be
\JJ\vert_{g.f.}=\eps_{-}+\sum_{a,b=1}^{N}\sum_{j=0}^{p-1}
W_{jj}^{ab} M^{jj}_{ab}\equiv \eps_{-}+\WW
\ee
where $W_{jj}^{ab}$ are the (unknown) generators of the gauge invariant 
 polynomials.
 
 In other words, there is a unique set of parameters $\lambda_{\mathbf x}$ 
 such that the gauge transformations (\ref{gaugtr}) leads 
 \JJc\ to \JJgf. These parameters are polynomials in the original 
 $J^{ab}_{jm}$, hence the generators $W_{jj}^{ab}$. Since they generate the gauge invariant 
 polynomials, the  $W_{jj}^{ab}$'s close (polynomially) under the PB: 
 they generate the $\cw(gl(Np),N.sl(p))$ algebra. The Lie algebra structure of 
 this \cw-algebra is given by the PB (\ref{PBgln}), together 
 with the knowledge of the polynomials $W_{jj}^{ab}$. Unfortunately, the 
 complete expression of these polynomials is difficult to obtain in 
 the general case, so that different technics have been developed to 
 compute the PB of the \cw-algebra, without knowing the exact 
 expression of the polynomials $W_{jj}^{ab}$. 
 
 There is essentially two different ways of defining the Poisson 
brackets of the $\cw_p(N)$ algebras: through 
the Dirac brackets, or using the so-called soldering procedure. We 
will need them both,  and describe them in the following.

\subsection{Dirac brackets\label{dirac}}
It can be shown that the first class constraints together with the 
gauge fixing form a set of second class constraints, \ie that 
if $\Phi=\{\phi_{\alpha}\}_{\alpha\in I}$ is the set of all constraints, we have
\be
\Delta_{\alpha\beta}=\{\phi_{\alpha},\phi_{\beta}\} \mbox{ is invertible: } 
\sum_{\gamma\in I}\Delta_{\alpha\gamma}\bar\Delta^{\gamma\beta}=
\delta_{\alpha}^\beta \mbox{ where } 
\bar\Delta^{\alpha\beta}\equiv(\Delta^{-1})_{\alpha\beta}
\ee
Together with a set of second class constraints occurs the notion of
 Dirac brackets which are constructed in such a way that they are 
compatible with these constraints:
\be
\{X,Y\}_{*}\sim \{X,Y\}-\sum_{\alpha,\beta\in 
I}\{X,\phi_{\alpha}\}\Delta^{\alpha\beta}\{\phi_{\beta},Y\}
\ \ \forall X,Y
\ee
where the symbol $\sim$ means that one has to apply the constraints on 
the right hand side {\em once the Poisson Brackets have been computed}.
The compatibility of the Dirac brackets with the constraints reflects 
in the following property
\be
\{X,\phi_{\alpha}\}_{*}\sim 0\ \ \forall \alpha\in I, \forall X
\ee

Then, the Poisson brackets of the \cw-algebra are defined as the Dirac 
brackets of the unconstrained generators $J_{jj}^{ab}$:
\be
\{W_{j}^{ab},W_{\ell}^{cd}\}\equiv 
\{J_{jj}^{ab},J_{\ell\ell}^{cd}\}_{*}
\ee

In the case we are considering, the matrix $\Delta$ 
take the form
\bea
\Delta_{jm;k\ell}^{ab;cd} &=& \{J_{jm}^{ab},J_{k\ell}^{cd}\}\ \ \ \forall 
a,b,c,d,j,k;\ \forall m<j;\ \forall\ell<k\\ 
 &=& (-1)^m \frac{j(j+1)-m(m+1)}{2} 
\frac{\eta_{j}}{\eta_{1}}\delta_{j,k}\delta_{m+\ell+1,0}\delta^{bc}\delta^{ad}+\\
&& + <j,m;k,\ell\vert t,t>\left(\delta^{bc} 
J^{ad}_{tt}-(-1)^{j+m+k+\ell}\delta^{ad}J^{cb}_{tt}\right)\\
&=& (-1)^m \frac{j(j+1)-m(m+1)}{2} 
\frac{\eta_{j}}{\eta_{1}}\delta_{j,r}\delta_{m+s+1,0}\delta^{be}\delta^{af}
\left( {\bf 1} - \widehat{\Delta}\right)_{k\ell;rs}^{cd;ef}  \label{Delta}\\
\widehat{\Delta}_{jm;k\ell}^{ab;cd} &=& 
\frac{2\eta_{1} <j,-m-1;k,\ell\vert 
t,t>}{\eta_{j}(j(j+1)-m(m+1))}\left(\rule{0ex}{1.2em}(-1)^{m}\,\delta^{ac} 
J^{bd}_{tt}+(-1)^{j+k+\ell}\delta^{bd}J^{ca}_{tt}\right)\ 
\eea
The form (\ref{Delta}) shows that $\Delta$ is invertible, for the 
matrix $\widehat{\Delta}$ is nilpotent: due to the Clebsch-Gordan 
coefficient $<j,-m-1;k,\ell\vert t,t>$, we have 
$(\widehat{\Delta})^{2p-1}=0$. Hence, we deduce 
\be
\bar\Delta^{jm;k\ell}_{ab;cd} = \label{invDelta}
(-1)^{m+1} \frac{j(j+1)-m(m+1)}{2} \,
\frac{\eta_{1}}{\eta_{j}}\ 
 \sum_{n=0}^{2p-1}(\widehat{\Delta}^n)_{j,-m-1;k\ell}^{ba;cd}  
\ee
where we have set 
\be
(\widehat{\Delta}^0)_{j,m;k\ell}^{ab;cd}=
\delta_{j,k}\delta_{m,\ell}\delta^{ac}\delta^{bd}
\ee
Once $(\Delta)^{-1}$ is known, one can compute the Dirac brackets. 
Unfortunately, in practice, (\ref{invDelta}) is difficult to achieve, 
and only partial results are obtained using the Dirac 
brackets. 

\subsection{Soldering procedure}
The calculation of the Poisson brackets of the \cw-algebra can be 
achieved through another way, called the soldering 
procedure\cite{ORaf}, see also \cite{RS} in the case of finite 
\cw-algebras.
It is not our aim to show the equivalence of this approach with the 
previous (Dirac) procedure. We give here jut a flavor of it in the 
context of $\cw_{p}(N)$-algebras. 

In the soldering procedure, the idea is to view the (adjoint) action 
of the $\cw_{p}(N)$ algebra on itself as a "residual" action of the whole 
$gl(Np)$ algebra on the currents, residual in the sense 
that is "respects" the constraints that have been imposed. In other 
words, among all the transformations induced by the (enveloping 
algebra of) $gl(Np)$, we look for the ones that do not affect the 
form $\JJ\vert_{g.f.}$: these will be the transformations induced by 
the $\cw_{p}(N)$-algebra. 

In the present paper, thanks to the basis explicited in the appendix 
\ref{notations}, we will be able to synthetically present (and solve) 
this procedure in the case of $\cw_{p}(N)$ algebras.

More precisely, the action of $gl(Np)$, with parameter 
$\lambda=\sum_{j,m;a,b}\lambda^{ab}_{jm}M^{jm}_{ab}$,
can be written
\be
\delta_{\lambda}\JJ=\{tr(\lambda \JJ),\JJ\}=[\lambda,\JJ ]
\ee
where $\{\ ,\ \}$ is the PB (on the $J$'s) and $[\ ,\ ]$ is the 
commutator (of $Np$x$Np$ matrices).
Within all these transformations, 
we look for the ones which preserve the form of $\JJ\vert_{g.f.}$:
\be
\delta_{\lambda}(\JJ\vert_{g.f.})=\delta_{\lambda}\WW=\sum_{j;a,b}(\delta_{\lambda} 
J_{jj}^{ab})M^{jj}_{ab}
\ee
This constrains the parameters $\lambda_{ab}^{jm}$, and 
only $N^2p$ of them are left free: they correspond to the parameters 
of the \cw-transformation. 

Explicitly, the calculation 
$[\lambda,\eps_{-}+\WW]=\delta_{\lambda}\WW$ leads to
\bea
\lambda^{j,m+1} &=& \sum_{k,r=0}^{p-1}\sum_{\ell=-k}^{k-1}\left(
\lambda^{k\ell}W_r <k,\ell;r,r\vert j,m> - W_r \lambda^{k\ell}
<r,r;k,\ell\vert j,m>\right)\nonumber\\
\mbox{for}&& -j\leq m< j \label{eq:lbd}\\
\delta_{\lambda} W_{j} &=&  \sum_{k,r=0}^{p-1}\sum_{\ell=-k}^{k-1}\left(
\lambda^{k\ell}W_r <k,\ell;r,r\vert j,j> - W_r \lambda^{k\ell}
<r,r;k,\ell\vert j,j>\right)\label{eq dW}
\eea
where $\lambda_{j,m}=\sum_{a,b}\lambda^{ab}_{jm}M^{jm}_{ab}$,
$W_{j}=\sum_{a,b} W_{j}^{ab} M^{jj}_{ab}$ and the products are 
matricial products.

The system (\ref{eq:lbd}) is strictly triangular in $\lambda_{jm}$ 
with respect to the gradation 
$gr(\lambda_{jm})=j+m$. Indeed, the Clebsch-Gordan coefficients ensure 
that $\vert j-r\vert\leq k\leq j+r$ and $\ell+r=m$, so that 
$gr(\lambda_{k,\ell})=k+\ell\leq j+m<j+m+1=gr(\lambda_{j,m+1})=$ in (\ref{eq:lbd}). 
Thus, all the $\lambda$'s are expressible in terms 
of the $\lambda_{j,-j}$ parameters.

\subsubsection{Calculation of $\{W_{0}^{ab},W_{j}^{cd}\}$}
As a start up, we consider the variation of $W_{0}$. In that case, 
one has only to look at (\ref{eq dW}), which reads:
\bea
\delta_{\lambda} W_{0} &=&  \sum_{k,r=0}^{p-1}\left(
\lambda^{k,-k}W_r <k,-k;r,r\vert 0,0> - W_r \lambda^{k,-k}
<r,r;k,-k\vert 0,0>\right)\\
&=& \sum_{k}^{p-1}(-1)^k\frac{\eta_{k}}{\eta_{0}}\, [\lambda^{k,-k},\,W_r]
\eea
Thus, we get the equation:
\be
\sum_{j}\widetilde\lambda_j\{W_{j},W_{0}\} = \frac{1}{p}\,\sum_{j}
[\widetilde\lambda_j,\,W_r]
\ee
where $\widetilde\lambda_j=(-1)^j \eta_{j}\lambda_{j,-j}$.
Hence, we are directly led to the PB:
\begin{equation}
    \{W_{0}^{ab},W_{j}^{cd}\}= \frac{1}{p}\,(\delta^{bc} W_{j}^{ad} -\delta^{ad} 
    W_{j}^{cb})
\end{equation}

\subsubsection{Calculation of $\{W_{1}^{ab},W_{j}^{cd}\}$}
Now, focusing on the variation of $W_1$ and using the results 
(\ref{res2:1}-\ref{res2:3}), we are led to
\bea
\delta_{\lambda}W_1 &=& 
\sum_{j}c_{j}\left(\frac{1}{j(2j-1)}\left[\lambda_{j-1,1-j},W_{j}\right] - 
\left[\lambda_{j,1-j},W_{j}\right]_{+}+\right.\nonumber\\
&&\left. -\frac{(j+1)(p-j-1)(p+j+1)}{2j+3}\left[\lambda_{j+1,1-j},W_{j}\right]\right)
\label{delW1}
\eea
where $c_{j}$ has been defined in (\ref{res2:1}) and $[\ ,\ ]$ (resp. $[\ ,\ ]_{+}$)
stands for the commutator (anti-commutator) of $Np$x$Np$ matrices.
Then, solving the equation (\ref{eq:lbd}) for $m=-j,1-j$, and plugging the 
result into (\ref{delW1}) gives
\bea
\sum_{j}\widetilde\lambda_j\{W_{j},W_{1}\} &=& {\textstyle \frac{3}{p(p^2-1)}}
\left(\ \sum_{j=1}^{p-1}\frac{j(p^2-j^2)}{2j+1} \right.
[\widetilde\lambda_{j-1},W_{j}]+
\sum_{j=1}^{p-1}\sum_{s=j}^{p-1}\left[\left[\widetilde\lambda_{s},
W_{s-j}\right],W_{j}\right]_{+}
\nonumber\\
 &&\hspace{10ex} +\sum_{j=0}^{p-1}\sum_{s=j+1}^{p-1}\frac{s-j-1}{2j+1} 
\left[\left[\widetilde\lambda_{s-1},W_{s-j-1}\right]_{+},W_{j}\right] +\\
&&\hspace{10ex} \left. - \sum_{j=0}^{p-1}\sum_{t=j+1}^{p-1}\sum_{s=t}^{p-1}\frac{1}{t(2j+1)} 
[[[\widetilde\lambda_{s},W_{s-t}],W_{t+1-j}],W_{j}]\right)\nonumber
\eea
In component, we get the following PB:
\be
\begin{array}{l}
\displaystyle
\{W_{1}^{ab},W_{j}^{cd}\} \ =\ {\frac{3}{p(p^2-1)}}
\left[\ \rule{0ex}{1.7em}\frac{(j+1)(p^2-(j+1)^2)}{2j+3}\ \left(\ \rule{0ex}{2.2ex}
\delta^{cb}W_{j+1}^{ad}-\delta^{ad}W_{j+1}^{cb}\ \right)\right.+\\[1.2em]
\displaystyle
\mbox{\hspace{1.2em}}+\ j\left(\rule{0ex}{2.2ex}
\delta^{cb}(W_{0}W_{j})^{ad}-\delta^{ad}(W_{j}W_{0})^{cb}
\ +\   W_{j}^{cb}W_{0}^{ad}-W_{j}^{ad}W_{0}^{cb}\right) 
+\\[1.em]
\displaystyle\mbox{\hspace{1.2em}}
+\sum_{s=1}^{j}\ (1+\frac{j-s}{2s+1})\ \left( \delta^{cb}(W_{s}W_{j-s})^{ad}
-\delta^{ad}(W_{j-s}W_{s})^{cb}\rule{0ex}{2.2ex}\ \right)+\\[1.42em]
\displaystyle\mbox{\hspace{1.2em}}
+\sum_{s=1}^{j}\ (1-\frac{j-s}{2s+1})\ \left(
W_{j-s}^{ad}W_{s}^{cb}-W_{s}^{ad}W_{j-s}^{cb}\rule{0ex}{2.2ex}\ \right) +\\[1.42em]
\displaystyle\mbox{\hspace{2.4ex}}
-\ \sum_{s=0}^{j-1}\sum_{t=s+1}^{j} \frac{1}{t(2s+1)}\left( 
\rule{0ex}{1.5em}
\delta^{cb}(W_{s}W_{t-s-1}W_{j-t})^{ad}- 
\delta^{ad}(W_{j-t}W_{t-s-1}W_{s})^{cb}\right.+\\[0.24em]
 \displaystyle
\mbox{\hspace{12em}}+ W_{j-t}^{ad}(W_{t-s-1}W_{s})^{cb}
- (W_{s}W_{t-s-1})^{ad}W_{j-t}^{cb}+\\[0.24em]
 \displaystyle\rule{0ex}{1.5em}
\mbox{\hspace{12em}}+ W_{t-s-1}^{ad}(W_{j-t}W_{s})^{cb}
-(W_{s}W_{j-t})^{ad}W_{t-s-1}^{cb}+\\[0.24em]
\displaystyle
 \mbox{\hspace{12em}}+ W_{s}^{ad}(W_{j-t}W_{t-s-1})^{cb}
-(W_{t-s-1}W_{j-t})^{ad}W_{s}^{cb}\left.\left.
\rule{0ex}{1.5em}\right)\rule{0ex}{1.7em}\right] 
\end{array}
\label{PB1j}
\ee

\sect{Comparison between truncated Yangians and finite 
\cw-algebras\label{reslt}}

We have seen that the  
truncated Yangians are a deformation of a truncated loop algebra based on $gl(N)$.
We show below that $\cw_{p}(N)$ is also a deformation of this 
algebra, and that these two deformations coincide. We use here 
the notions presented in appendix \ref{coho}.

We work at the classical (Poisson brackets) level.

\subsection{$\cw_p(N)$ as a deformation of a truncated 
loop algebra}
To see that the $\cw_p(N)$ is a deformation of a truncated 
loop algebra based on $gl(N)$, we modify the constraints to 
\be
\JJ=\frac{1}{\hbar} \eps_{-}+\sum_{a,b=1}^{N}\sum_{j=0}^{p-1}\sum_{0\leq m\leq j} 
J_{jm}^{ab} M^{jm}_{ab}
\ee
These constraints are equivalent to the previous ones as soon as 
$\hbar\neq0$. With these new constraints, the matrix $\Delta$ and its inverse 
read
\bea
(\Delta_{\hbar})_{jm;k\ell}^{ab;cd} &=&  \frac{1}{\hbar}\,(-1)^m \frac{j(j+1)-m(m+1)}{2} 
\frac{\eta_{j}}{\eta_{1}}\delta_{j,k}\delta_{m+\ell+1,0}\delta^{bc}\delta^{ad}
\left( 1 - \hbar\widehat{\Delta}_{k\ell;rs}^{cd;ef} \right) \nonu
(\bar\Delta_{\hbar})^{jm;k\ell}_{ab;cd} &=& \hbar\ (-1)^{m+1} \frac{j(j+1)-m(m+1)}{2} 
\frac{\eta_{1}}{\eta_{j}}\
 \sum_{n=0}^{2p-1} \hbar^n\, (\widehat{\Delta}^n)_{j,-m-1;k\ell}^{ba;cd} 
\eea
Then, computing the Dirac brackets associated to these 
new constraints, one finds
\bea
\{J_{jj}^{ab},J_{\ell\ell}^{cd}\}_{\hbar}&=&\{J_{jj}^{ab},J_{\ell\ell}^{cd}\}-
\sum_{efgh;kmrs}\{J_{jj}^{ab},J_{km}^{ef}\}(\bar\Delta_{\hbar})^{km;rs}_{ef;gh}
\{J_{rs}^{gh},J_{\ell\ell}^{cd}\}\\
&=& 
\delta^{bc}J_{j+\ell,j+\ell}^{ad}-\delta^{ad}J_{j+\ell,j+\ell}^{cb}-\hbar 
P_{\hbar}(J)
\eea
where $P_{\hbar}(J)$ (polynomial in the $J_{jj}^{ab}$ which is computed using 
$\bar\Delta_{\hbar}$ as in section \ref{dirac}) has only positive (or null) 
powers of $\hbar$. 
This clearly shows that the $\cw_p(N)$ algebra 
is a deformation of the algebra generated by $W_{j}^{ab}\equiv 
J_{jj}^{ab}$ and with defining (undeformed) Poisson brackets:
\bea
\{W_{j}^{ab},W_{\ell}^{cd}\}_{0} &=& \delta^{bc}W_{j+\ell}^{ad}-\delta^{ad}W_{j+\ell}^{cb}
\ \mbox{ if }j+\ell<p \\
&=& 0 \ \mbox{ if }j+\ell\geq p
\eea
One recognizes in this algebra a (enveloping) loop algebra based on 
$gl(N)$ quotiented by the relations $W_{j}^{ab}=0$ if $j\geq p$. In 
other words, this algebra is nothing but a truncated loop algebra, and 
the \cw-algebra is a deformation of it.

\subsection{Identification of $\cw_p(N)$ and $\ynp{p}$}

We have already seen that the truncated Yangians as well as 
the \cw-algebras we consider are both a deformation of a truncated 
loop algebra:
\bea
\{W_j^{ab},W_\ell^{cd}\}_{1} &=& \{W_{j}^{ab},W_{\ell}^{cd}\}_{0} + 
\sum_{n=1}^\infty\hbar^n\vph^{W}_n(W_j^{ab},W_\ell^{cd}) \ \ 0\leq 
j,\ell\leq p-1\\
\{\bar T_m^{ij},\bar T_n^{kl}\}_{2} &=& \{\bar T_{m}^{ij},\bar T_{n}^{kl}\}_{0} + 
\sum_{r=1}^\infty\hbar^r\vph^{T}_r(\bar T_{m}^{ij},\bar T_{n}^{kl})\ \ 0\leq 
m,n\leq p-1 
\eea
where the cochains $\vph_{n}^W$ and $\vph_{r}^T$ obey 
(\ref{delta:phi1}-\ref{delta:phin}). 
The undeformed PB $\{\ ,\ \}_{0}$ are identical (via the 
identification\footnote{The shift $j\rightarrow j-1$ in the 
identification is due to a difference of convention between 
\cw-algebras and Yangians: in the former case, the index $j$ denotes 
the underlying $sl(2)$ representation, while in the latter $j$ is the 
exponent of $u$ in the formal series (\ref{eq:defT}).}  
$W_{j}^{ab}\equiv \bar T_{j-1}^{ab}$) and correspond to 
the truncated loop algebra. Thus, we have two deformed PB $\{\,,\,\}_1$
and $\{\,,\,\}_2$, and all we need is to show that the cochains
 $\vph_{n}^W$ and $\vph_{n}^T$ coincide $\forall n$. To prove 
that it is indeed the case, we need the following properties:

{\lem \label{lemCocycle}Let $gl(N)_{p}$ be the loop algebra based on 
$gl(N)$, truncated at order $p$, and $u_j^{ab}$ $(j<p)$ the 
corresponding generators. A 2-cocycle $\vph$ with values in 
$\cu(gl(N)_{p})$ is completely determined once one knows 
$\vph(u_0^{ab},u_j^{cd})$ and $\vph(u_1^{ab},u_j^{cd})$, 
$\forall a,b,c,d=1,\ldots,N$ and $\forall 
j=0,\ldots,p-1$}

\null

\prf
We prove this lemma recursively.
We write the cocycle condition for a triplet 
$(u_1^{A},u_j^{B},u_k^{C})$,  using indices 
$A,B,C=1,\ldots,N^2$ in the adjoint representation, and 
the commutation relations of $gl(N)_{p}$:
\bea
&&{f^{AB}}_{D}\, \vph(u_{1+j}^{D},u_{k}^{C})+{f^{BC}}_{D} \, 
\vph(u_{k+j}^{D},u_{1}^{A})+{f^{CA}}_{D}\, \vph(u_{1+k}^{D},u_{j}^{B})
=\nonumber\\
&&\ =\ \{u_{1}^{A}, \vph(u_{j}^{B},u_{k}^{C}) \}+\{u_{j}^{B}, 
\vph(u_{k}^{C},u_{1}^{A}) \}+\{u_{k}^{C}, \vph(u_{1}^{A},u_{j}^{B}) \}
\eea
It can be rewritten as
\bea
\gamma_{2}\ \vph(u_{1+j}^{A},u_{k}^{B}) &= &
{f_{A}}^{CD} {f_{DB}}^{E}\, \vph(u_{1}^{C},u_{k+j}^{E})+
{f_{A}}^{CD} {f_{DB}}^{E}\, \vph(u_{j}^{C},u_{k+1}^{E})+\\
&&+{f_{A}}^{CD}\left(\{u_{k}^{B}, \vph(u_{j}^{C},u_{1}^{D}) \}+
\{u_{j}^{C}, \vph(u_{1}^{D},u_{k}^{B}) \}+
\{u_{1}^{D}, \vph(u_{k}^{B},u_{j}^{C}) \}\right)\nonumber
\label{eq:recur}
\eea
where $\gamma_{2}\neq 0$ is the value of the second Casimir operator 
in the adjoint representation.

For $j=1$, (\ref{eq:recur}) allows to compute 
$\vph(u_{2}^{D},u_{k}^{C})$ $\forall C,D$ and $\forall k\geq 1$ once 
$\vph(u_{1}^{D},u_{k}^{C})$ $\forall C,D$ and $\forall k$ is known. 

Suppose now 
that we know $\vph(u_{j}^{A},u_{k}^{B})$ for $1\leq j<\ell_{0}$ and 
$\forall k$. Then, 
(\ref{eq:recur}) for $j=\ell_{0}-1$ allows to compute 
$\vph(u_{\ell_{0}}^{A},u_{k}^{B})$  $\forall k$. 

Thus, apart from the values $\vph(u_{0}^{A},u_{k}^{B})$ we are able 
to compute all the expressions $\vph(u_{j}^{A},u_{k}^{B})$. This ends 
the proof.
\finprf
\begin{prop} \label{propBas}
    There exist two sets of generators $\{{}^{\pm}\bar W_{j}^{ab}\}_{j=0,\ldots}$ in 
$\cw_{p}(N)$ such that 
\bea
&&\{{}^{\pm}\bar W_1^{ab}\,,\,{}^{\pm}\bar W_{j}^{cd}\}= \delta^{cb}\,
{}^{\pm}\bar W^{ad}_{j+1}-\delta^{ad}\, {}^{\pm}\bar W^{cb}_{j+1} +
\bar W^{cb}_{0}\, {}^{\pm}\bar W^{ad}_{j}- {}^{\pm}\bar W^{cb}_{j}\, \bar W^{ad}_{0}\nonu
&&\forall a,b,c,d=1,\ldots,N\ ;\  \forall j\geq 1
\label{eq.PBn}\\
&&\{\bar W_0^{ab}\,,\,{}^{\pm}\bar W_{j}^{cd}\}= \delta^{cb}\, 
{}^{\pm}\bar W^{ad}_{j}-\delta^{ad}\, {}^{\pm}\bar W^{cb}_{j} \nonumber
\eea
The generators $^{\pm}\bar W_{j}^{ab}$ are polynomial of degree 
$(j+1)$ in 
the original ones $W_j^{ab}$ and are recursively defined by
\bea
{}^{\pm}\bar W_{j,\pm}^{ab} &=&\sum_{n=1}^{j+1}{}^{\pm}\bar 
W_{j,(n)}^{ab}\ =\ \sum_{n=1}^{j+1}\sum_{\mid \vec{s}\mid =j+1-n} 
{}^{\pm}\alpha_{\vec{s}}^{n,j}\, (W_{s_{1}}\ldots W_{s_{n}})^{ab}\ 
1<n\mbox{ and } 1<j\nonu
{}^{\pm}\bar W_{1}^{ab} &=& \pm\frac{p(p^2-1)}{6}\ W_1^{ab}+\frac{p(p\pm1)}{2}\ 
(W_{0}W_{0})^{ab}
\label{eq.Wjn}\\
\bar W_{0}^{ab} &\equiv& {}^{+}\bar W_{0}^{ab}\ =\ {}^{-}\bar W_{0}^{ab}
\ =\ p\ W_0^{ab} 
\nonumber
\eea
for some numbers ${}^{\pm}\alpha_{\vec{s}}^{n,j}$. The summation on $\vec{s}$ 
is understood as a summation on $n$ 
positive (or null) integers $(s_{1},\ldots,s_{n})\equiv\vec{s}$ such that 
$\vert \vec{s}\vert=\sum_{i=1}^m s_{i}=j+1-n$. 

The subsets $\{{}^{\pm}\bar W_{j}^{ab}\}_{j=0,\ldots,p-1}$ form two 
bases of $\cw_{p}(N)$, the other generators $\{{}^{\pm}\bar 
W_{j}^{ab}\}_{j\geq p}$ been polynomials in the basis elements.
\end{prop}

\null

\prf
We first remark that the form (\ref{eq.Wjn}) clearly shows that the 
$p$ first generators are independent, and thus form a basis. The other 
ones must then be polynomials in any basis.

We prove the relation (\ref{eq.PBn}) by a recursion on $j$. 
It is easy to compute 
that the definitions  are such that (\ref{eq.PBn}) 
is satisfied for $j=1$. For the recursion, we fix a basis $^{+}\bar W_{j}^{ab}$ or 
$^{-}\bar W_{j}^{ab}$  (the proof is obviously independent of the 
choice), and write it $\bar W_{j}^{ab}$.

We suppose that we have found generators $\bar W_{j}^{ab}$ for $j\leq j_{0}$ 
such that (\ref{eq.PBn}) is satisfied. This implies that we have:
\[
N\, \bar W^{cd}_{j_{0}+1} = 
\{\bar W_1^{ca},\bar W_{j_{0}}^{ad}\}-
\bar W^{aa}_{0}\,\bar W^{cd}_{j_{0}}+\bar W^{cd}_0
\, \bar W^{aa}_{j_{0}}+\delta_{{cd}}\bar W^{aa}_{j_{0}+1}
\]
where we have used implicit summation on repeated $gl(N)$ indices.
Then, we get
\bea
N \{\bar W_1^{ab},\bar W^{cd}_{j_{0}+1}\} &= & 
\{\bar W_1^{ab},\{\bar W_1^{ce},\bar W_{j_{0}}^{ed}\}\}-
\{\bar W_1^{ab},\bar W^{ee}_{0}\,\bar W^{cd}_{j_{0}}-
\bar W^{cd}_{0}\, \bar W^{ee}_{j_{0}}\}
+\delta_{{cd}}\{\bar W_1^{ab},\bar W^{aa}_{j_{0}+1}\}\nonu
&=& \{\{\bar W_1^{ab},\bar W_1^{ce}\},\bar W_{j_{0}}^{ed}\}
+\{\bar W_1^{ce},\{\bar W_1^{ab},\bar W_{j_{0}}^{ed}\}\}
+\delta_{{cd}}\{\bar W_1^{ab},\bar W^{aa}_{j_{0}+1}\}+\nonu
&& -\bar W^{ee}_{0}\,
\{\bar W_1^{ab},\bar W^{cd}_{j_{0}}\}+
\bar W^{cd}_{0}\, \{\bar W_1^{ab},\bar W^{ee}_{j_{0}}\}+
\bar W^{ee}_{j_{0}}\, \{\bar W_1^{ab},\bar W^{cd}_{0}\}\nonu
&=& \{\bar W_1^{cb},\bar W^{ad}_{j_{0}+1}\}-
\{\bar W_2^{cb},\bar W^{ad}_{j_{0}}\}+\
\bar W^{ab}_{1}\,\bar W^{cd}_{j_{0}}-
\bar W^{cd}_{1}\, \bar W^{ab}_{j_{0}}\nonu
&&+N\left(\bar W^{cb}_{0}\,\bar W^{ad}_{j_{0}+1}-
\bar W^{ad}_{0}\, \bar W^{cb}_{j_{0}+1}\right) +
\delta^{cb}\, A^{ad}-\delta^{ad}\, B^{cb}+\delta^{cd}\, C^{ab}\nonumber
\eea
with the notation
\beano
A^{ad} &=& \{\bar W_2^{ae}, \bar W^{ed}_{j_{0}}\}+
\bar W^{ad}_{0}\,\bar W^{ee}_{j_{0}+1}-
\bar W^{ee}_{0}\, \bar W^{ad}_{j_{0}+1}+
\bar W^{ad}_{1}\,\bar W^{ee}_{j_{0}}-
\bar W^{ee}_{1}\, \bar W^{ad}_{j_{0}}\\
B^{cb} &=& \{\bar W_1^{ce}, \bar W^{eb}_{j_{0}+1}\}+
\bar W^{cb}_{0}\,\bar W^{ee}_{j_{0}+1}-
\bar W^{ee}_{0}\, \bar W^{cb}_{j_{0}+1}\\
C^{ab} &=& \{\bar W^{ab}_{1}\, , \,\bar W^{ee}_{j_{0}}\}+
 [\bar W_{0}, \bar W_{j_{0}+1}]^{ab}
\eeano
It remains to compute $\{\bar W_2^{cb},\bar W^{ad}_{j_{0}}\}$. 
This is done using the same technics as above:
\[
N\, \bar W_2^{cb} = \{\bar W_1^{ce},\bar W_1^{eb}\}-
\bar W^{ee}_{0}\,\bar W^{cb}_{1}+\bar W^{cb}_{0}\, \bar W^{ee}_{1}+
\delta^{cb}\, \bar W_2^{ee}
\]
so that we have
\bea
N\{\bar W_2^{cb},\bar W^{ad}_{j_{0}}\} &=& 
-\left(\rule{0ex}{1.2em}\{\bar W_1^{ab},\bar W^{cd}_{j_{0}+1}\}+
\{\bar W_1^{cd},\bar W^{ab}_{j_{0}+1}\}\right)
+N\left(\rule{0ex}{1.2em}\bar W^{ab}_{1}\,\bar W^{cd}_{j_{0}}-
\bar W^{cd}_{1}\, \bar W^{ab}_{j_{0}}\right)+
\nonu
&& +N\left(\rule{0ex}{1.2em}
\bar W^{ab}_{0}\,\bar W^{cd}_{j_{0}+1}-
\bar W^{cd}_{0}\, \bar W^{ab}_{j_{0}+1}\right)+\delta^{ab}B^{cd}+
\label{W2Wj}
\\
&&+\delta^{cd}\left(\rule{0ex}{1.2em}-\{\bar W_{j_{0}+1}^{ae}\, , \,\bar W^{eb}_1\}+
\bar W^{ee}_{0}\,\bar W^{ab}_{j_{0}+1}-
\bar W^{ab}_{0}\, \bar W^{ee}_{j_{0}+1}\right)+\nonu
&&+\delta^{cb}\left(\rule{0ex}{1.2em}\{\bar W^{ee}_{2}\, , \,\bar W^{ad}_{j_{0}}\}-
 [\bar W_{0}\, , \, \bar W_{j_{0}+1}]^{ad}-
 [\bar W_{1}\, , \, \bar W_{j_{0}}]^{ad}\right)
 \nonumber
\eea
Then, a recurrent use of these two brackets leads to the result:
\bea
\{\bar W_1^{ab}\, , \,\bar W^{cd}_{j_{0}+1}\} &=&
\bar W^{cb}_{0}\,\bar W^{ad}_{j_{0}+1}-
\bar W^{ad}_{0}\, \bar W^{cb}_{j_{0}+1} +
\delta^{cb}\ \bar W^{ad}_{j_{0}+2}-\delta^{ad}\ \bar 
W^{cb}_{j_{0}+2}+\nonu
&&-
\frac{\delta^{ad}\delta^{cb}}{N(N^2-1)}\left(\{\bar W_1^{ee}\,,\,\bar 
W^{ff}_{j_{0}+1}\}-N\{\bar W^{ef}_{1}\, ,\,\bar 
W^{fe}_{j_{0}+1}\}\right)+\nonu
&&+\frac{\delta^{cd}}{N}\left(\rule{0ex}{1.5em}\{\bar W_1^{ab}\,,\,\bar 
W^{ee}_{j_{0}+1}\}+[\bar W_{0}\,,\,\bar W_{j_{0}+1}]^{ab}+\right.
\label{PB.fin}\\
&&\left.\phantom{+\delta^{cd}(} +\frac{\delta^{ab}}{N(N^2-1)}
\left(\rule{0ex}{1.2em}\{\bar W_1^{ee}\, , \,\bar 
W^{ff}_{j_{0}+1}\}-N\{\bar W^{ef}_{1}\, , \, \bar 
W^{fe}_{j_{0}+1}\}\right)\right) 
\nonumber
\eea
for some polynomials $\bar W^{ad}_{j_{0}+2}$.

Finally, we remark that the forms (\ref{eq.Wjn}) and the PB (\ref{PB1j}) 
clearly show that the PB
$\{\bar W_{{1}}^{ab},\bar W^{cd}_{j}\}$ does not contain terms 
proportional to $\delta^{ab}$ or $\delta^{cd}$. Moreover, a direct 
calculation, using (\ref{PB1j}), shows that 
\[
\{\bar W_{{1}}^{aa}, P^{bb}_{j}\}=0\ ,\ \ \forall\ 
P^{cd}_{j}(W)=\sum_{n=1}^{j+1}\sum_{\mid \vec{s}\mid =j+1-n} 
\beta_{\vec{s}}^{n,j} (W_{s_{1}}\ldots W_{s_{n}})^{cd}
\]
This is enough to show that 
the two last lines in the PB (\ref{PB.fin}) identically vanish.

Hence, we can deduce that the PB must be of the form
\be
\{\bar W_1^{ab},\bar W^{cd}_{j_{0}+1}\} =
\bar W^{cb}_{0}\,\bar W^{ad}_{j_{0}+1}-
\bar W^{ad}_{0}\, \bar W^{cb}_{j_{0}+1} +
\delta^{cb}\ \bar W^{ad}_{j_{0}+2}-\delta^{ad}\ \bar W^{cb}_{j_{0}+2}
\ee
which is exactly (\ref{eq.PBn}), so that the recursion on $j$ is  
proven. 
\finprf
We have computed the first and last terms ($\forall j\geq0$) 
that appear in the definition 
(\ref{eq.Wjn}): 
\bea
{}^{\pm}\bar W_{j,(1)}^{ab} &=& (\pm1)^j\, (j!)^2\, {p+j\choose 
2j+1}\ W_j^{ab} \\
{}^{-}\bar W_{j,(j+1)}^{ab} &=& {p\choose 
j+1}\ (\underbrace{W_0\cdots W_{0}}_{j+1})^{ab} 
\label{Wmj0}\\
{}^{+}\bar W_{j,(j+1)}^{ab} &=& {p+j\choose 
j+1}\ (\underbrace{W_0\cdots W_{0}}_{j+1})^{ab}
\label{Wpj0}
\eea
{\coro The change of generators between $\{{}^{+}\bar W_j^{ab}\}$ and 
$\{{}^{-}\bar W_j^{ab}\}$ is given by:
\be
{}^{\pm}\bar W_j^{ab} = \sum_{n=1}^{j} (-1)^{j+n+1}\sum_{\mid \vec{s}\mid =j+1-n}
({}^{\mp}\bar W_{s_{1}}\ldots {}^{\mp}\bar W_{s_{n}})^{ab}
\label{eq.chgB}
\ee
\label{chgB}}
\prf
Using the expression (\ref{eq.Wjn}) for $j=1$ and the PB (\ref{eq.PBn}),
one computes that 
\be
\left\{{}^{\pm}\bar W_{1}^{ab}\, ,\, {}^{\mp}\bar W_{j}^{cd} \right\}=
\delta^{bc}\left(\rule{0ex}{1.2em}(\bar W_{0} {}^{\mp}\bar W_{j})^{ad}-{}^{\mp}\bar 
W_{j+1}^{ad}\right) - \delta^{ad} \left(\rule{0ex}{1.2em}({}^{\mp}\bar W_{j}\bar 
W_{0})^{cb}-{}^{\mp}\bar W_{j+1}^{cb}\right)
\ee
Then, a direct calculation shows that indeed the expression (\ref{eq.chgB}) 
satisfies (\ref{eq.PBn}).
\finprf
{\coro The basis $\{{}^-\bar W_j^{ab}\}$ is such that ${}^-\bar W_j^{ab}=0$ for $j\geq 
p$. In the basis $\{{}^+\bar W_j^{ab}\}$, all the ${}^+\bar W_j^{ab}$ 
generators ($j\geq p$) are not vanishing.}
 
\null

\prf From the PB (\ref{eq.PBn}) it is clear that it is sufficient to 
show that $\bar W_p^{ab}=0$.
Writing this PB for $j=p$ and using the form
\be
\bar W_p^{ab}=\sum_{n=2}^{p+1}\sum_{\mid \vec{s}\mid =p+1-n} 
\alpha_{\vec{s}}^{n,p} (\bar W_{s_{1}}\ldots \bar W_{s_{n}})^{ab}
\ee
one gets only two possibilities for the $\alpha$'s:
\be
\alpha_{\vec{s}}^{n,p}= (-1)^n A \ \mbox{ with }\ A=0\mbox{ or }1
\ee
If $A=0$, then $\bar W_p^{ab}=0$ while if $A=1$,
the change of basis given in the corollary \ref{chgB} shows that 
in the other basis we have $\bar W_p^{ab}=0$. Hence, we have to 
determine which basis corresponds to $\bar W_p^{ab}=0$.

Looking at the expressions (\ref{Wpj0}) and (\ref{Wmj0}), 
one concludes that ${}^-\bar W_p^{ab}=0$ while ${}^+\bar 
W_j^{ab}\neq0,\ \forall\, j$.
\finprf
In the following, we choose for $\cw_{p}(N)$ the $\{{}^-\bar W_j^{ab}\}$ 
basis and omit the superscript $-$ for the generators.
Now, from above, it is easy to show:
\newpage
{\theor \label{W-Y} The $\cw_{p}(N)$ algebra is the truncated 
Yangian $\ynp{p}$.}

\null

\prf
Let us first remark that the two algebras have identical (in fact 
undeformed) PB on the 
couples $(\bar W_0^{ab},\bar W_j^{cd})$, which proves that the 
cochains $\vph_{n}^W$ and 
$\vph_{n}^T$ coincide (in fact vanish) on these points.

Moreover, the property \ref{propBas} shows that the cochains $\vph_{n}^W$ and 
$\vph_{n}^T$ coincide on the couples $(\bar W_1^{ab},\bar W_j^{cd})$. 
Since $\vph_{1}$ is a cocycle, this is enough (using lemma 
\ref{lemCocycle}) to prove 
that $\vph_{1}^T$ and $\vph_{1}^W$ are identical. 

Now, suppose that we have proven that $\vph_{n}^W$ and 
$\vph_{n}^T$ are identical for $n<n_{0}$. Then, eq. 
(\ref{delta:phin}) fixes $\vph_{n_{0}}^W$ and 
$\vph_{n_{0}}^T$,  up to a cocycle:
\beano
\vph_{n_{0}}^W &=& \vph_{n_{0}}+\xi_{n_{0}}^W\\
\vph_{n_{0}}^T&=& \vph_{n_{0}}+\xi_{n_{0}}^T
\eeano
where $\vph_{n_{0}}$ is a function of the cochains
$\vph_{n}^W=\vph_{n}^T$, $n<n_{0}$. 
But property \ref{propBas} shows that the two cocycles $\xi_{n_{0}}^W$ and 
$\xi_{n_{0}}^T$ coincide on the couples $(\bar W_1^{ab},\bar 
W_j^{cd})$, which proves that they are identical (due to lemma 
\ref{lemCocycle}). 

Thus, $\vph_{n_{0}}^W$ and 
$\vph_{n_{0}}^T$ are identical, and we have proven recursively the property.
\finprf

\subsubsection{Quantization}
We have shown that truncated Yangians and \cw-algebras coincide at the 
classical level. It remains to show that it is still true at the 
quantum level. Fortunately, an algebra morphism between Yangians and 
\cw-algebras has already been given in \cite{RS}, at classical and quantum 
levels. This relation was not sufficient to establish the 
identification
between \cw-algebras and truncated Yangians, since all the 
 horizontal arrows involved in the diagram 
 \be
 \begin{array}{ccc}
     \displaystyle Y(N) & \longrightarrow & \cw_{p}(N) \\[1.2ex]
      \updownarrow & & \updownarrow ? \\[1.2ex]
    \displaystyle Y(N) & \longrightarrow & Y_{p}(N) 
\end{array}\label{isoWY}
\ee
 are not isomorphisms. Hence the calculations done in this paper. 
 
 However, once the relation (between $Y_{p}(N)$ and $\cw_{p}(N)$) has been 
 established 
 at the classical level, we can use the result of \cite{RS} to promote 
 it at the quantum level. More precisely, now that we can identify 
 the $\cw_{p}(N)$ algebra with $Y_{p}(N)$ at the classical level, we 
 can use the results of \cite{RS} at the quantum level: it has been 
 established that any quantization of $\cw_{p}(N)$ still obey to the 
 Drinfeld relation, and hence the homomorphism still exists at the 
 quantum level.
 
 Thus, theorem \ref{W-Y} is valid both at classical and quantum 
 level, and the figure \ref{isoWY} is correct (without question mark).

 Let us remark that in the proof we have establish, we have 
 constructed \cw-algebras as deformations of a truncated loop algebra 
 and identified them with the 
truncated Yangians,\ie truncations of deformed loop algebras. 
Denoting by $\cl(gl(N))$ the loop algebra defined on $gl(N)$, and by 
$\cl(gl(N))_{p}$ its truncation, the above sentence can be pictured as 
the following commutative diagram:
\be
\begin{array}{ccccc}
    & & Y(N) & & \\
    & \nearrow_{\hbar} & & \searrow^{p} &\\
    \cl(gl(N)) & & & & Y_{p}(N) \equiv \cw_{p}(N)\\
    & \searrow^{p} & & \nearrow_{\hbar} &\\
     & &\cl(gl(N))_{p}& & 
 \end{array}
 \ee
 where $\nearrow_{\hbar}$ stands for a deformation, and $\searrow^{p}$ for a 
 truncation (at level $p$).
   
\section{Applications\label{appl}}
\subsection{$R$-matrix for \cw-algebras}
The above construction allows us to associate the \cw-algebras to 
the $R$-matrix of the Yangian, the difference between these two 
algebras lying in the modes development of $T(u)$: in both cases, the 
development is done in powers of $u^{-1}$, but for the Yangian 
it is an infinite series, while the development is truncated to a 
polynomial for the \cw-algebra. Explicitly, the presentation of the 
$\cw_{p}(N)$-algebra take the form:
\[
   R(u-v)T_{1}(u)T_{2}(v) = T_{2}(v)T_{1}(u)R(u-v)
\ \mbox{ with }\ \left\{\begin{array}{l}\displaystyle
T(u)=1+\sum_{n=1}^p\sum_{a,b=1}^{N^2} u^{-n} E_{ab}\, T^{ab}_{n}
\\[1.3em]
\displaystyle R(x)=1\otimes 1-\frac{1}{x}P_{12} \end{array}\right.
\]

Let us remark that this procedure is similar to the "factorization 
procedure" which leads from the elliptic algebra $\ca_{q,p}(N)$ to 
the Sklyanin algebra $S_{q,p}(N)$ \cite{skly} (see also \cite{clad} 
for more examples about factorizations). In all cases, one chooses for $T(u)$ a special 
dependence in $u$ to get a finite algebra: this special dependence is 
nothing but a coset by some of the modes of $T(u)$. In all the 
examples, the Hopf structure of the starting algebra does not 
survive to this quotient.

Note  that the $R$-matrix presentation of the \cw-algebras 
provides an {\em exhaustive} set of commutation relations among the 
$\cw_{p}(N)$ generators for generic $N$ and $p$, while, up to now, 
a complete set of commutation relation was known only for 
a small number of \cw-algebras. 

Let us also remark that the $R$-matrix presentation allows to define 
the \cw-algebras without any reference to the underlying $gl(Np)$ 
algebra, and thus is a more "abstract" definition.

\subsection{Irreducible representations of $\cw_{p}(N)$-algebras}
Once again, the $R$-matrix presentation provides a very natural 
framework for the classification of \cw-representations\footnote{We 
thank P. Sorba for drawing our attention to this point.}. It is based 
on the notion of evaluation representations, as it appears in the 
Yangian context (see section \ref{sec.eval}).
In fact, this classification was done in \cite{Chered}, in the context of 
(truncated) Yangians. We have the following theorem:
\begin{theor}\label{thm.rep}
    {\bf Finite dimensional irreducible representations of 
    $\cw_p(N)$}\\    
    Any finite dimensional irreducible representation of the $\cw[gl(Np),N.gl(p)]$ 
    algebra is isomorphic to an evaluation representation or to the 
    subquotient of tensor product of  at most $p$ evaluation representations.
\end{theor}
\prf

By evaluation representations for $\cw_{p}(N)$ algebra, we mean the 
definitions \ref{def.eval} and \ref{tens.eval} with the change 
$T^{ab}_{r}\rightarrow\, W^{ab}_{r-1}$ (\ie the evaluation 
representations of the truncated Yangian). The property (\ref{tensW}) 
clearly shows that the (subquotient of) 
tensor product of $n$ evaluation representations is 
a representation of the truncated Yangian as soon as $n\leq p$. It 
also shows that if it is irreducible for the Yangian, then it is also 
irreducible for the truncated Yangian and that they are finite 
dimensional.

Now conversely, an irreducible representation $\pi$ of the $\cw_{p}(N)$ 
algebra can be lifted to a representation of the whole Yangian by 
setting $\pi(T^{ij}_{(r)})=0$ for $r>n$. It is then obviously irreducible 
for the Yangian, and thus is isomorphic to the tensor product of 
evaluation representations.
\finprf
We remark that the theorem \ref{thm.rep} allows to construct any 
(finite dimensional) representation of $\cw_{p}(N)$ in term of $p$ 
representations of $gl(N)$ (including trivial representations). 
This is exactly what one obtains from the 
so-called "Miura transformation" that appears in the context of 
\cw-algebras. Indeed, this transformation allows to construct a 
representation of the $\cw(\cg,\ch)$-algebra using representations of
$\cg_{0}$, the zero-grade 
subalgebra  of $\cg$. In the case of $\cw_{p}(N)$, we get 
$\cg_{0}=N.gl(p)$, and hence 
need $N$ representations of $gl(p)$, as it is stated in theorem \ref{thm.rep}.

Finally, as for Yangians, we have the following characterization 
(proved using above theorem and the characterization for Yangians):

\begin{coro}
    The irreducible finite-dimensional representations of $\cw_{p}(N)$ 
    are in one-to-one correspondence with the families 
    $\{P_{1}(u),\ldots,P_{N-1}(u),\rho(u)\}$ where $P_{i}$ are monic 
polynomials of degree $m_{i}$ such that $\sum_{i}m_{i}\leq p$, 
and $\rho(u)=1+\sum_{n>0} d_{n}u^{-n}$.
\end{coro}

\subsection{Generalization to $\synp{r}{p}$ truncated Yangians.}
As well as we have defined truncated Yangians based on $Y(N)$, the same 
construction can be done for each of the $\synp{r}{}$ Hopf algebras, to construct
$\synp{r}{p}$ algebras ($r\leq p$): these truncated Yangians will correspond to the 
quotient of $\cw_{p}(N)$ by $\cd_{r}$, which is a part of the 
$\cw_{p}(N)$-center (see below).  

In particular, the
$\cw(sl(Np), N.sl(p))$ algebra usually encountered in the literature 
is nothing but the truncation $\synp{1}{p}$. For this algebra, one 
sees that there exist two algebra homomorphisms: $Y(gl(N))\rightarrow 
\cw(sl(Np), N.sl(p))$ and $Y(sl(N))\rightarrow \cw(sl(Np), N.sl(p))$. 
The second one corresponds to the case given in \cite{RS}.

More generally, we have the following property:
\begin{prop}
    There is an algebra homomorphism from $\synp{r}{p+q}$ to 
    $\synp{r+s}{p}$, for any values of $p,q,r,s=0,1,\ldots,\infty$.
\end{prop}
\prf
It is a trivial composition of algebra and Hopf algebra 
homomorphisms, as it is visualized in figure \ref{fig}.
\finprf

\begin{figure}[htp]
$$
\begin{array}{cccccccccc}
  \begin{array}{c} Y(N)\\ \vertequiv \\ \synp{0}{\infty} \end{array} & 
  \!\Rightarrow\! & \synp{1}{\infty}  & \cdots & \Rightarrow & 
  \synp{p}{\infty} & \cdots & \Rightarrow & 
  \begin{array}{c} \synp{\infty}{\infty}\\ \vertequiv \\ Y(sl(N))\end{array}  
 \\
  \downarrow & \rule{0ex}{1.2em} & \downarrow &&& \downarrow &&&
  \d2dots\hspace{4.1em}
 \\
 \vdots & & \vdots &&& \vdots &&\d2dots& 
 \\
 \downarrow & & \downarrow &&& \downarrow &\d2dots&&
 \\
 \begin{array}{c} \ynp{p}\\ \vertequiv \\ \cw_{p}(N)\end{array} &  
 \!\Rightarrow\! & 
 \begin{array}{c} \synp{1}{p}\\ \vertequiv \\ 
     \mbox{\small{$\cw[sl(Np),N.sl(p)]$}}\end{array}  
 & \cdots & \Rightarrow & \synp{p}{p} &&&\\
 \rule{0ex}{1.2em}\downarrow & & \downarrow &&& 
 \d2dots\hspace{4.1em} &&& 
 \\
 \vdots & & \vdots &&\d2dots&&&& 
 \\
 \downarrow & & \downarrow &\d2dots&&&&&
 \\
 \begin{array}{c} \ynp{1}\\ \vertequiv \\  \cu(gl(N))\end{array} & \!\Rightarrow\! & 
 \begin{array}{c}\synp{1}{1} \\ \vertequiv \\ \cu(sl(N))\end{array} &&&&&&
 \\
 & & &&&&&&
 \\
\downarrow &\d2dots& &&&&&&
 \end{array}
$$
$\ynp{0}\equiv\{1\}$
\caption{\it The vertical links $\downarrow$ correspond to the 
truncations (algebra homomorphisms),  
while the horizontal links $\Rightarrow$ are associated to coset by central elements
(Hopf algebra homomorphisms). }\label{fig}
\end{figure}
\newpage
\subsubsection{Finite dimensional irreducible representations of 
$\synp{r}{p}$ algebras}
Starting from the theorem \ref{thm.rep} and using cosets by central 
elements, it is easy to get
\begin{coro} Any finite-dimensional irreducible representation of 
the $\synp{r}{p}$ algebra is obtained from the subquotient of tensor product of at 
most $p$ evaluation representations, quotiented by $r$ 
constraints on the generators of $\cd_{r}$.

The finite-dimensional irreducible representations of 
the $\synp{r}{p}$ algebra are in one-to-one correspondence with the 
families     $\{P_{1}(u),\ldots,P_{N-1}(u),\rho(u)\}$ where $P_{i}$ are monic 
polynomials of degree $m_{i}$ such that $\sum_{i}m_{i}\leq p$, 
and $\rho(u)=1+\sum_{n=0}^r d_{n}u^{-n}$.
\end{coro}
In particular, in the case of the $\cw(sl(Np),N.sl(p))$ algebra, we 
obtain the result given in \cite{RS} for $N=2$.
\subsection{Center of $\cw_{p}(N)$ algebras.}
From the definition of $\cw(gl(Np),N.sl(p))$ algebras, one already 
knows that their center contains the Casimir operators of $gl(Np)$, 
since, being central, these operators are obviously gauge invariant.
Hence the dimension of the center is at least $Np$. However, it was 
not proved (to our knowledge) that its dimension is exactly $Np$. 
Fortunately, 
the center of the truncated Yangians $Y_{p}(N)$ has been determined in 
\cite{Chered}:

\null

{\bf{Property:}} {\em 
    A basis of the $Y_{p}(N)$ center is given by all
the coefficients of the principal part of the following generating 
function
\be
H(x)=\sum_{w\in S_{N}}\sum_{i=1}^N \sum_{r_{i}=0}^{p-1}\ (-1)^{sg(w)}\ 
T^{w(1)1}_{r_{1}}\, T^{w(2)2}_{r_{2}}\, \cdots\, T^{w(N)N}_{r_{N}}
\ \prod_{j=1}^N\left( \frac{(x-j)^{p-1-r_{j}}}{\prod_{k=1}^p 
(x-j-u_{k})}\right)
\ee
where $S_{N}$ is the symmetric group and $T^{ab}_{r}$ are the Yangian 
generators.}

\null

Looking at the poles of $H(x)$, it is easy to see that there are 
exactly $Np$ poles (including multiplicities). A basis for this 
center (using quantum determinant) was also given in \cite{yangTr}.
Hence, using this property and the above remark,
we can deduce
\begin{coro}
    The center of $\cw_{p}(N)$ is $Np$-dimensional and is given by 
    $\cd_{Np}/\ci_{p}$. A basis of this center is canonically associated to 
    the Casimir operators of $gl(Np)$.
 \end{coro}
 Let us remark that the $p$ first Casimir operators can be chosen as 
 elements of the $\cw_{p}(N)$ basis, while the next $p(N-1)$ ones are 
 polynomials in the basis generators. 
 Note that a different way 
 to get these central generators has been given in \cite{Wcom}. It 
 uses a determinant formula for $gl(Np)$ expressed for $\JJ_{gf}$, 
 namely:
 \be
 \mbox{det}(\JJ_{gf}-\lambda\, \II)=(-1)^{Np}\lambda^{Np}\, +\, \sum_{n=0}^{Np-1} 
 C_{Np-n}\lambda^n
 \ee

More generally, the same reasoning leads to the following center for 
$\synp{r}{p}$: 
\be
Z(\synp{r}{p})=\cu(d_{r+1},\ldots,d_{pN})/\ct_{p}
\ee
It is generated by the 
last $(Np-r)$ independent Casimirs of $gl(Np)$.

\sect{Conclusion\label{conclu}}
We have shown that the finite $\cw(gl(Np),N.sl(p))$ algebras are 
nothing but truncated Yangians $Y(gl(N))_{p}$, \ie
coset of the Yangian $Y(gl(N))$ by the relations $T^{ab}_{(n)}=0$ 
for $n\geq p$. The resulting coset is an algebra, but the Yangian Hopf 
structure does not survive to the quotient.
This property enlightens 
the algebra homomorphism between Yangians and finite 
\cw-algebras, and which was given in \cite{RS}.
Using this property, we 
have been able to present these \cw-algebras as exchange algebras, 
with the help of the Yangian R-matrix. This more abstract 
presentation is not linked to an Hamiltonian reduction, as were 
usually defined the \cw-algebras. It could be of some help in the 
seek of a geometrical interpretation of \cw-algebras.
As a consequence, we have also 
given a complete classification of the finite dimensional irreducible 
representations for these \cw-algebras. This classification completes
the one given in \cite{RS} for $\cw(sl(2n),2.sl(n))$ algebras. 
Physically, one can hope to construct lattice models associated to 
\cw-algebras, starting from  models with Yangian symmetry. 

Now that the relation between Yangians and \cw-algebras is 
well-understood, one can hope to construct R-matrices for general 
\cw-algebras: work is in progress in this direction.
Conversely, one can think of generalizing the notion of Yangian as 
certain limits of $\cw(\cg,\ch)$ algebras in which a (quasi) Hopf structure 
can be recovered. This would provide a wide class of new types of 
quantum groups.

Let us also remark that two other approaches for Yangians and \cw-algebras 
could be related. On the one hand, one can construct Yangians as 
the projective limit of the centralizer of $gl(n)$ in $\cu(gl(m+n))$ 
\cite{Ol} (see also \cite{OM}), and on the other hand, some finite 
\cw-algebras (of type $\cw(gl(2n),n.sl(2))$) have been realized as 
commutants of a $gl(2n)$ parabolic subalgebra in a certain 
localization of $\cu(gl(2n))$ \cite{Wcom}.  It seems to us quite natural 
to look for a global description of these two point of view. 

Of course, the case of conformal \cw-algebras (\ie extensions of the 
Virasoro algebra) has to be considered. It could be related to a
multi-parametric generalization of Yangians.  Would such a 
generalization be possible, one could think of an "RTT 
presentation" of Virasoro algebra: this would allow to relate "usual" 
\cw-algebras with the deformed 
\cw-algebras presentation, a link which is not clear up to now, since 
two different deformed algebras can be constructed \cite{rims,lapth}.
Note finally that the construction of some conformal \cw-algebras 
(such as the Virasoro and the $\cw_{3}$ algebras) as commutant in a 
localization of an affine Kac-Moody algebra (see above paragraph) as 
been already achieve \cite{Wcom}: this could be a way to generalize the 
notion of Yangians, using the centralizer construction.

\null

{\bf\Large Acknowledgments}

\null

We warmly thank Daniel Arnaudon, Michel Bauer and Paul Sorba for 
fruitful and clarifying discussions and Alexander Molev 
for his enlightening comments about representations of Yangians.

\appendix
\section{General settings on $gl(Np)$\label{notations}}
We have gathered here the notations and properties we  need about
$gl(Np)$ algebra.

We consider the $gl(Np)$ algebra in its fundamental representation 
($Np\,$x$\, Np$ matrices), and take a basis adapted to the 
decomposition with respect to the $sl(2)$ algebra principal in 
$N.sl(p)\equiv \underbrace{sl(p)\oplus\ldots\oplus sl(p)}_{N}$. 
This decomposition makes naturally appear the "factorization" 
$gl(Np)=gl(N)\otimes gl(p)$, valid in the fundamental representation, 
and the $sl(2)$ principal in $sl(p)$. 

\subsection{The principal embedding of $sl(2)$ in $sl(p)$}
We will denote by  $M_{j,m}$ (with $-j\leq m\leq j$ and $1\leq j\leq 
p$) the $p\,$x$\, p$ matrices resulting from the 
decomposition in $sl(2)$ multiplets:
\begin{eqnarray}
    {[e_{+},M_{j,m}]} &=& \frac{j(j+1)-m(m+1)}{2} M_{j,m+1} \label{e+Mjm}\\
    {[e_{-},M_{j,m}]} &=&  M_{j,m-1} \label{e-Mjm}\\
    {[e_{0},M_{j,m}]} &=& m\ M_{j,m} \label{e0Mjm}\\
    {[e_{0},e_{\pm}]} &=& \pm\ e_{\pm} \ \mbox{ and }\ 
    {[e_{+},e_{-}]}\ =\ e_{0}
\end{eqnarray}
where $e_{\pm,0}$ are the generators of the $sl(2)$ algebra principal 
in $sl(p)$. The normalizations in (\ref{e+Mjm}-\ref{e-Mjm}), although 
not symmetric, are adapted to the \cw-algebra framework. When working 
with $gl(p)$ instead of $sl(p)$, we will add the $j=0$ generator, 
proportional to the identity matrix.

The decomposition of $M_{j,m}$ in terms of the $p\,$x$\, p$ matrices $E_{ab}$ reads
\be
\begin{array}{lll}
    \displaystyle
    M_{j,m} = \sum_{k=1}^{p-m} a_{j,m}^k E_{k,k+m} &\mbox{ with }& 
    \displaystyle a_{j,m}^k = \sum_{i=0}^{j-m} (-1)^{i+j+m} {j-m\choose i} 
    a_{j,j}^{k-i}\\[1.4em] 
   \displaystyle
     \mbox{ for } 0\leq m\leq j &&
\end{array}
\ee
\be
\begin{array}{lll}
    \displaystyle
     M_{j,m} = \sum_{k=1}^{p+m} a_{j,m}^k E_{k-m,k} &\mbox{ with }&
    \displaystyle a_{j,m}^k = \sum_{i=0}^{j-m} (-1)^{i+j+m} {j-m\choose i} 
    a_{j,j}^{k-i-m}\\[1.4em]  
   \displaystyle
    \mbox{ for }-j\leq m\leq0 &&\\[1.4em] 
\end{array}
\ee
\be
    a_{j,j}^{k} = \frac{(k+j-1)! (p-k)!}{(k-1)! (p-k-j)!} 
\ee
The generators $e_{\pm,0}$ of the $sl(2)$ algebra are proportional to the $M_{1,m}$ 
generators:
\bea
    e_{+} &=& \sum_{k=1}^{p-1} \frac{k(p-k)}{2} E_{k,k+1} \ = \
    \half M_{1,1} \nonu
    e_{0} &=& \sum_{k=1}^{p} (\frac{p+1}{2}-k) E_{k,k} \ =\  - \half 
    M_{1,0}\\
    e_{-} &=& \sum_{k=1}^{p-1} E_{k+1,k} \ = \ - \half M_{1,-1} \nonumber
\eea
Let us remark that we have the following generating function for the 
coefficients $^{(p)}a_{j,m}^k$ (where $(p)$ refers to the $gl(p)$ 
algebra under consideration):
\bea
&&^{(p)}a_{j,m}^k=\frac{j!}{p!\, (k-1)!\, (j-m)!}\left[\frac{d^p}{du^p}
\frac{d^j}{dz^j}\frac{d^{j-m}}{dy^{j-m}}\frac{d^{k-1}}{dx^{k-1}}
\ a(x,y,z;u)\right]_{\begin{array}{l}
\scriptstyle x=y=0\\[-1.2ex] \scriptstyle z=u=0\end{array}}\\
&&\mbox{ with }\ 
a(x,y,z,u)=\frac{u}{\left(\rule{0ex}{1.em}1+y(1-x)\right)
\left(\rule{0ex}{2.4ex}1-u[1+z+x(1-u)]\right)}
\eea
The scalar product is given by
\begin{eqnarray}
    \eta_{j,m;\ell,n} &=& (M_{j,m}, M_{\ell,n}) =tr(M_{j,m}\cdot M_{\ell,n}) = 
    (-1)^m \eta_{j}\ 
    \delta_{j,\ell}\delta_{m+n,0} \\
    &&\mbox{ with } \eta_{j}= (2j)! (j!)^2 
    {p+j\choose 2j+1}
\end{eqnarray}
where the dot stands for the matrix product and $tr$ is the trace of 
matrices (in the $p$-dimensional representation). 

In the following, we will need the Clebsch-Gordan like coefficients 
given by:
\begin{equation}
    M_{j,m}\cdot M_{\ell,n} = \sum_{r=\vert j-\ell\vert}^{j+\ell}\,\sum_{s=-r}^r 
    <j,m;\ell,n\vert r,s> M_{r,s} \label{defClb}
\end{equation}
As for usual Clebsch-Gordan 
coefficients, one can prove (using commutators by $e_{\pm,0}$) 
that $r$ must be in $[\vert 
j-\ell\vert\, , \, j+\ell]$ and that
$s$ must be equal to $-m-n$. However, 
since we are in the fundamental of $sl(p)$, the coefficients will be 
truncated in such a way that only the values $r\leq p$ are kept in 
the decomposition (\ref{defClb}). We 
will still call them Clebsch-Gordan coefficients.

Using the scalar product, one can compute these 
coefficients to be
\begin{equation}
    <j,m;\ell,n\vert r,s> = \frac{(-1)^s}{\eta_{r}} tr\left( M_{j,m}\cdot 
    M_{\ell,n}\cdot M_{r,-s}\right)
\end{equation}

\subsection{Few results about the Clebsch-Gordan like coefficients}
Using the cyclicity of the trace, one shows 
\begin{equation}
    <j,m;\ell,n\vert r,s> = (-1)^{s+m}\frac{\eta_{j}}{\eta_{r}} 
    <\ell,n;r,-s\vert j,-m>= (-1)^{s+n}\frac{\eta_{\ell}}{\eta_{r}} 
    <r,-s;j,m\vert \ell,-n>
\end{equation}

We will also use the property
\begin{equation}
    <j,m;\ell,n\vert r,s> = \frac{(j-m)! (\ell-n)!(r+s)!}{(j+m)! (\ell+n)!(r-s)!} 
    <\ell,-n;j,-m\vert r,-s>
\end{equation}
where the coefficients are due to the non-symmetric basis we have 
chosen.

With these two properties, one can compute:
\bea
<r,-r;k,k\vert j,-j>  &=&  (-1)^k\ \frac{\eta_{j+k}}{\eta_{j}} \ 
\delta_{r,j+k}\label{res1:1} \\ 
<k,k;r,-r\vert j,-j> &=& (-1)^k\ \frac{\eta_{j+k}}{\eta_{j}} \ 
\delta_{r,j+k}\\
<r,1-r;k,k\vert j,1-j>  &=&   (-1)^k\ \frac{j}{j+k}\frac{\eta_{j+k}}{\eta_{j}} 
\ \delta_{r,j+k} \\ 
<k,k;r,1-r\vert j,1-j>  &=&  (-1)^k\ \frac{j}{j+k}\frac{\eta_{j+k}}{\eta_{j}} 
\ \delta_{r,j+k}\\
<r,-r;k,k\vert j,1-j>  &=&  (-1)^{k+1}\ k\, j\, \frac{\eta_{j+k-1}}{\eta_{j}} 
\ \delta_{r+1,j+k} \\ 
<k,k;r,-r\vert j,1-j>  &=& (-1)^k\ k\, j\, \frac{\eta_{j+k-1}}{\eta_{j}} 
\ \delta_{r+1,j+k}\label{res1:6}
\eea
We will also need  the following coefficients:
\bea
<k,k;k,1-k\vert 1,1>  &=&  (-1)^{k+1}\ \frac{\eta_{k}}{\eta_{1}}\ 
\equiv\ c_{k} 
\label{res2:1}\\  
<k,1-k;k,k\vert 1,1>  &=&  -c_{k} \\  
<k,k;k-1,1-k\vert 1,1>  &=&  \frac{1}{k(2k-1)}\ c_{k}
\\ 
<k-1,1-k;k,k\vert 1,1>  &=&  \frac{1}{k(2k-1)}\ c_{k}
\\ 
<k+1,1-k;k,k\vert 1,1>  &=&  -\frac{(k+1)(p^2-(k+1)^2)}{2k+3}\ c_{k}
\\  <k,k;k+1,1-k\vert 1,1>  &=&  -\frac{(k+1)(p^2-(k+1)^2)}{2k+3}\ 
c_{k} 
\label{res2:3}
\eea

\subsection{Basis for $gl(Np)$}
We can use the above basis of $gl(p)$ to construct a basis for 
$gl(Np)$. Using the  $N\,$x$\, N$ matrices $E_{ab}$, the generators 
$\Upsilon_{ab}^{jm}$ of $gl(Np)$ in the fundamental will be represented by
\begin{equation}
    \pi_{F}(\Upsilon_{ab}^{jm})=M_{ab}^{jm}=E_{ab}\otimes M^{jm}
\end{equation}
The generators of the $sl(2)$ algebra principal in $N.sl(p)$ are then
\begin{equation}
    \eps_{\pm,0} = 1_{N}\otimes e_{\pm,0}
 \end{equation}
 where $e_{\pm,0}$ are the $p\,$x$\, p$ matrices defined above. 
 We have the following commutation relations
 \begin{eqnarray}
    {[\eps_{+},M_{ab}^{j,m}]} &=& \half(j(j+1)-m(m+1)) M_{ab}^{j,m+1} \\
    {[\eps_{-},M_{ab}^{j,m}]} &=&  M_{ab}^{j,m-1} \\
    {[\eps_{0},M_{ab}^{j,m}]} &=& m M_{ab}^{j,m} \\
    {[\eps_{0},\eps_{\pm}]} &=& \pm \eps_{\pm} \ \mbox{ and }\ 
    {[\eps_{+},\eps_{-}]}\ =\ \eps_{0}
\end{eqnarray}
together with 
\begin{equation}
    {[M_{ab}^{00},M_{cd}^{00}]} \ =\ 
    \delta_{bc}M_{ad}^{00}-\delta_{ad}M_{cb}^{00}
\end{equation}
 This last commutator reveals the $gl(N)$ algebra which commutes with 
 the $sl(2)$ subalgebra under consideration.
 
More generally, the product law (in the fundamental representation) reads
\begin{equation}
    M_{ab}^{j,m}\cdot M_{cd}^{\ell,n} = 
    \delta_{bc}\sum_{r=\vert j-\ell\vert}^{ j+\ell}\sum_{s=-r}^r 
    <j,m;\ell,n\vert r,s> M_{ad}^{r,s}
\end{equation}
which leads to the following commutation relations (valid in the 
abstract algebra):
\begin{equation}
    [\Upsilon_{ab}^{j,m}, \Upsilon_{cd}^{\ell,n}] = 
    \sum_{r=\vert j-\ell\vert}^{ j+\ell}\sum_{s=-r}^r \left(
    \rule{0ex}{2.4ex}
    \delta_{bc} <j,m;\ell,n\vert r,s> \Upsilon_{ad}^{r,s} - 
    \delta_{ad} <\ell,n;j,m\vert r,s> \Upsilon_{cb}^{r,s} \right)
\end{equation}
The scalar product is
\begin{equation}
    \eta_{ab,cd}^{j,m;\ell,n} = (\Upsilon_{ab}^{j,m}, \Upsilon_{cd}^{p,n}) 
    =tr(M_{ab}^{j,m}\cdot M_{cd}^{p,n}) = 
    \delta_{a,d}\ \delta_{b,c}\ \eta^{j,m;\ell,n}
\end{equation}

\section{Deformations and cohomology\label{coho}}
 We include here  some definitions (in the context 
of Chevalley cohomology) to be self-content. For more details about 
deformations and their relation to cohomology, we refer to 
\cite{hoch} and ref. therein.

\subsection{Few words about Chevalley cohomology\label{cheval}}
We begin with an algebra \ca, and first introduce the space $C_{n}(\ca,\ca)$ 
of $n$-cochains with values in \ca, \ie 
skew-symmetric linear maps from $\wedge^n\ca$ to $\ca$.
The Chevalley derivation $\delta$ maps $n$-cochains to 
$(n+1)$-cochains as:
\beano
(\delta\chi_{n})(u_{0},\ldots,u_{n}) &=& \sum_{i=0}^n\ (-1)^i\ 
\left\{\rule{0ex}{2.1ex}u_{i}\ ,\ \chi_{n}(u_{0},u_{1},\ldots,\wh{u_{i}},\ldots,u_{n}) \right\}+ \\
&& + \sum_{0\leq i<j\leq n}\ (-1)^{i+j}\ 
\chi_{n}\left( \{u_{i}, u_{j}\},u_{0},u_{1},\ldots,\wh{u_{i}},\ldots,\wh{u_{j}},\ldots,u_{n}
\right)
\eeano
where, as usual, $\wh{u_{i}}$ means that $u_{i}$ has to be discarded 
in the list (or product, or sum, or whatever) we consider.

It can be shown that $\delta$ squares to zero:
\be
(\delta(\delta\chi_{n}))(u_{-1},u_{0},u_{1},\ldots,u_{n}) =0\ \ 
\forall u_{-1},u_{0},u_{1},\ldots,u_{n}\ ;\ \forall \chi_{n}\ ;\ 
\forall n
\ee
Thus, we introduce the cohomology associated to $\delta$, \ie we 
focus on Ker$\delta$. Elements of Ker$\delta$ are called cocycles, 
and we will see that they play a direct role in the deformation of 
Lie algebras. The space of $n$-cocycles (with values in \ca) is denoted 
$Z_n(\ca,\ca)$: Ker$\delta=\oplus_{n}Z_n(\ca,\ca)$.
Since $\delta^2=0$, we have 
Im$\delta\subset$Ker$\delta$: each $n$-cochain provides a $(n+1)$-cocycle. 
The elements $\delta\chi_{n}$  
correspond to "trivial" cocycles: they are called coboundaries, and 
the corresponding space denoted $B_n(\ca,\ca)$. The cohomology describes the 
non-trivial cocycles, \ie it is the space 
$H_n(\ca,\ca)=Z_n(\ca,\ca)/B_n(\ca,\ca)$, $H(\ca,\ca)=\oplus_{n} 
H_n(\ca,\ca)=$Ker$\delta/$Im$\delta$.
Due to its definition, the Chevalley cohomology is naturally 
associated to Lie algebras. When the cochains take values in $\CC$
instead of \ca, the space $H_{2}(\ca,\CC)$ classifies the non-trivial 
central extensions of \ca: see for instance \cite{KM} where central 
extensions of generalized loop algebras are classified and computed.
In the case we are considering, 
$C(\ca,\ca)$ is related to deformations of \ca.

\subsection{Deformations}
We start again with an algebra, with generators $u_{\alpha}$ 
($\alpha\in \Gamma$).
$$
\{u_{\alpha}, u_{\beta}\} = {f^{\alpha\beta}}_{\gamma} u_{\gamma}
$$
Actually, we will consider its enveloping algebra $\ca$, and introduce 
a deformation of it
$$
\{u_{\alpha}, u_{\beta}\}_{\hbar} = {f^{\alpha\beta}}_{\gamma} u_{\gamma}
+\sum_{n=1}^\infty \hbar^n \varphi_{n}(u_{\alpha}, u_{\beta})
$$
where the antisymmetric bilinear forms $\varphi_{n}$ take values in $\ca$:
they are all elements of $C_{2}(\ca,\ca)$.

Asking the bracket $\{.,.\}_{\hbar}$ to obey the Jacobi identity leads 
to the following equations:
\bea
\delta\vph_{1} &=& 0 \label{delta:phi1}\\
\hspace{-2.1em}\delta\vph_{n} &=& \hspace{-1.2ex}\sum_{j+k=n}\!\!\left(\rule{0ex}{2.4ex}\, \vph_{j}(\vph_{k}(u,v),w)+ 
\vph_{j}(\vph_{k}(v,w),u)+\vph_{j}(\vph_{k}(w,u),v) \right) 
\label{delta:phin}
\mbox{ for }  n>1 
\eea
where the operation $\delta$ defined in section \ref{cheval} has naturally appeared.

These equations indicate that $\vph_{1}$ is a cocycle, while 
$\vph_{n}$ is determined by the $\vph_{p}$'s ($p<n$) up to a cocycle.
Note that $\delta\vph_{n}$ is a coboundary, so that the $\vph_{p}$, 
$p<n$, must be such that the r.h.s. of (\ref{delta:phin}) is also a 
coboundary (it can be proven that this r.h.s. is indeed a cocycle, 
\ie is annihilated by $\delta$). If the third cohomological space is not 
trivial, the r.h.s. of (\ref{delta:phin}) may be a cocycle while being \underline{not} a 
coboundary: this leads to the usual assumption that the third cohomological 
space classify the obstructions to deformations. In other words, it 
could appear that, in the attempt to construct a deformation, the 
chosen $\vph_{p}, p<n$ are such that the l.h.s. of (\ref{delta:phin}) 
is a non-trivial cocycle, so that one cannot solve this equation at 
level $n$. In that case, the deformation would be ill-defined.

Fortunately, in the case 
we will consider below, we already know that we have well-defined 
deformations, and we have not to deal with a possible obstruction.

Note also that if $\vph_{n}$ is a coboundary 
$$
\vph_{n}(u_{\alpha},u_{\beta})=\delta\chi_{n}(u_{\alpha},u_{\beta})=
\{u_{\alpha},\chi_{n}(u_{\beta})\}-\{u_{\beta},\chi_{n}(u_{\alpha})\}-
\chi_{n}(\{u_{\alpha},u_{\beta}\})
$$
we can perform a change of basis 
$$\tilde{u}_{\alpha}=u_{\alpha}-\hbar^n \chi_{n}(u_{\alpha})
$$ 
such that in this new basis, the term in $\hbar^n$ has disappeared:
$$
\{\tilde{u}_{\alpha},\tilde{u}_{\beta}\}_{\hbar} = {f^{\alpha\beta}}_{\gamma}
\tilde{u}_{\gamma}+
\sum_{m=1}^{n-1}\hbar^m \varphi_{m}(\tilde{u}_{\alpha}, \tilde{u}_{\beta})
+\sum_{m=n+1}^{\infty} \hbar^m \tilde\varphi_{m}(\tilde{u}_{\alpha}, \tilde{u}_{\beta})
$$
where $\tilde\varphi_{m}, m>n$ are new cochains resulting from the 
change of variables. In that sense, a coboundary leads to a trivial deformation. 
However, one has to be careful that to "trivialize" the full 
deformation,  the change of basis has to be done recursively and the coboundarity 
at level $n$ has to be checked once the change of basis at 
level $n-1$ has been done (since the cochains are modified at higher 
order).

\end{document}